\documentclass[lang = american,12pt,oneside]{amsart} 

\usepackage{amsthm , amsmath, enumitem}
\usepackage{geometry , caption , graphicx}
\usepackage{hyperref , fontenc , inputenc}
\usepackage{array , babel , booktabs}
\usepackage{cite , float , footmisc}
\usepackage{iftex , indentfirst , kvoptions}
\usepackage{newtxmath , newtxtext , pdf14 , pdftexcmds}
\usepackage{ragged2e , url , xcolor , xpatch}

\usepackage{stmaryrd}
\usepackage{xypic}

\usepackage{environ}

\theoremstyle{plain}
\newtheorem*{theorem*}{Theorem}
\newtheorem*{proposition*}{Proposition}
\newtheorem*{lemma*}{Lemma}
\newtheorem*{conjecture*}{Conjecture}

\theoremstyle{definition}
\newtheorem*{example*}{Example}
\newtheorem*{remark*}{Remark}

\DeclareFontFamily{U} {MnSymbolA}{}

\DeclareFontShape{U}{MnSymbolA}{m}{n}{
  <-6> MnSymbolA5
  <6-7> MnSymbolA6
  <7-8> MnSymbolA7
  <8-9> MnSymbolA8
  <9-10> MnSymbolA9
  <10-12> MnSymbolA10
  <12-> MnSymbolA12}{}
\DeclareFontShape{U}{MnSymbolA}{b}{n}{
  <-6> MnSymbolA-Bold5
  <6-7> MnSymbolA-Bold6
  <7-8> MnSymbolA-Bold7
  <8-9> MnSymbolA-Bold8
  <9-10> MnSymbolA-Bold9
  <10-12> MnSymbolA-Bold10
  <12-> MnSymbolA-Bold12}{}

\DeclareSymbolFont{MnSyA} {U} {MnSymbolA}{m}{n}

\DeclareMathSymbol{\dashedleftarrow}{\mathrel}{MnSyA}{98}
\DeclareMathSymbol{\dashedrightarrow}{\mathrel}{MnSyA}{96}

\newcommand{\diag}{{\operatorname{diag}}}
\newcommand{\Gm}{\mathbb G_m}
\newcommand{\Ga}{\mathbb G_a}
\newcommand{\GL}{\operatorname{GL}}
\newcommand{\SL}{\operatorname{SL}}
\newcommand{\PGL}{\operatorname{PGL}}
\newcommand{\Sp}{\operatorname{Sp}}
\newcommand{\Spin}{\operatorname{Spin}}
\newcommand{\SO}{\operatorname{SO}}
\newcommand{\LG}{{^LG}}
\newcommand\aut{{\operatorname{aut}}}
\newcommand{\RTF}{\operatorname{RTF}}
\newcommand{\hor}{{\operatorname{hor}}}
\newcommand{\ver}{{\operatorname{ver}}}
\newcommand{\st}{{\operatorname{st}}}
\newcommand{\shear}{{\mathbin{\mkern-6mu\fatslash}}}

\newcommand{\withproofs}[1]{ #1}
\newcommand{\withoutproofs}[1]{} 

\numberwithin{equation}{section}

\begin{document}


\title{Spherical varieties, functoriality, and quantization}

\author{Yiannis Sakellaridis}

\address{Johns Hopkins University, 3400 N.\ Charles St., Baltimore, MD 21218, USA} \email{sakellar@jhu.edu}


\begin{abstract}
We discuss generalizations of the Langlands program, from  \emph{reductive groups} to the local and automorphic spectra of \emph{spherical varieties}, and to more general representations arising as ``quantizations'' of suitable Hamiltonian spaces. To a spherical $G$-variety $X$, one associates a \emph{dual group} ${^LG_X}$ and an \emph{$L$-value} (encoded in a representation of ${^LG_X}$), which conjecturally describe the local and automorphic spectra of the variety. This sets up a problem of functoriality, for any morphism ${^LG_X}\to {^LG_Y}$ of dual groups. We review, and generalize, Langlands' ``beyond endoscopy'' approach to this problem. Then, we describe the cotangent bundles of quotient stacks of the relative trace formula, and show that transfer operators of functoriality between relative trace formulas in rank 1 can be interpreted as a change of ``geometric quantization'' for these cotangent stacks.

\end{abstract}

\maketitle


\section{Integral representations of $L$-functions}

\subsection{Classical periods}

\subsubsection{} In his legendary 1859 paper \cite{Riemann}, Riemann proved the functional equation of the zeta function by representing it as the Mellin transform of a theta series
\[\pi^{-\frac{s}{2}} \Gamma\Big(\frac{s}{2}\Big) \zeta(s) = \int_0^\infty y^{\frac{s}{2}} \sum_{n=1}^\infty e^{-n^2\pi y} \, d^\times y.\]
The proof used the functional equation of the latter with respect to the substitution $y \leftrightarrow y^{-1}$, previously established by Jacobi and based on the Poisson summation formula. 

About 90 years later, Iwasawa, in a short announcement \cite{Iwasawa}, and Tate, in his thesis \cite{Tate}, reformulated this integral in the language of the adeles. The new formulation could be directly applied to the generalizations of the zeta function to arbitrary Dirichlet characters (by Dirichlet), or number fields (by Dedekind) and Gr\"ossencharacters (by Hecke), and  
clarified the meaning of the Euler factors of the zeta function, as Mellin transforms of Schwartz functions on the $p$-adic completions of $\mathbb Q$. Namely, we have an identity
\[ \zeta_p(s) = \int_{\mathbb Q_p^\times} \Phi_p(x) |x|^s \, d^\times x,
\]
where, for finite primes $p$, $\zeta_p(s) = (1-p^{-s})^{-1}$ and $\Phi_p = 1_{\mathbb Z_p}$, the characteristic function of the $p$-adic integers, is what we will call the \emph{basic Schwartz function} on $\mathbb Q_p$; the same interpretation extends to the ``Archimedean factor'' $\zeta_\infty (s) =  \pi^{-\frac{s}{2}} \Gamma\Big(\frac{s}{2}\Big)$ of the functional equation, with $\mathbb Q_\infty = \mathbb R$,  and $\Phi_\infty$ the Gaussian $e^{-\pi x^2}$.

\subsubsection{} Meanwhile, in 1936--37, Hecke \cite{Hecke1,Hecke2,Hecke3} had introduced what is today called the $L$-function of a modular form, generalized to non-holomorphic automorphic forms by Maass in 1944 \cite{Maass}. Recast in the adelic language by Jacquet and Langlands in their seminal 1970 work \cite{JL}, these $L$-functions, with appropriate Archimedean factors, can be represented as Mellin transforms 
\[ \int_{k^\times \backslash \mathbb A^\times} f \begin{pmatrix}
                                                     a \\ & 1 
                                                    \end{pmatrix} |a|^s \, d^\times a,\]
where $k$ denotes a number field, $\mathbb A$ its ring of adeles, and $f$ is a cuspidal automorphic form on $\GL_2(k)\backslash \GL_2(\mathbb A)$.

Shortly after Hecke, Rankin \cite{Rankin} and Selberg \cite{Selberg} discovered an integral representation for the $L$-function that carries their names, which today is seen as a special case of a Langlands $L$-function, attached to the tensor product representation 
\[ \check G=\GL_2 \times \GL_2 \xrightarrow{\otimes} \GL_4\]
of the Langlands dual of the group $G=\GL_2 \times \GL_2$. This integral is, on the surface, very different from the Mellin transforms of Riemann and Hecke, as it involves a pair of cusp forms and an Eisenstein series:
\[ \int_{\GL_2(\mathbb A)} f_1(g) f_2(g) E(g,s) \,  dg.\]

\subsection{The theta series of spherical varieties}

\subsubsection{} \label{thetaseries} The aforementioned works, and their adelic reformulations, led to an explosion of research around $L$-functions from the 70s onward, with numerous new integral representations discovered by Godement, Jacquet, Rallis, Piatetski-Shapiro, Gelbart, Shalika,  Waldspurger, Ginzburg, Bump, Friedberg, Garrett, and others \cite{Bump-RS}, combining elements from all of the methods above, such as the theta series (from Riemann), the ``period integrals'' over subgroups (from Hecke), and the Eisenstein series (from Rankin and Selberg). 

A uniform approach to many of these methods was proposed in \cite{SaRS}; it relies on the following ingredients:
\begin{itemize}[leftmargin=2em]
 \item A (suitable) \emph{affine spherical variety} $X$ for a group $G$ over a number field $k$; that is, $X$ is a normal, affine $G$-variety, with a dense orbit for the Borel subgroup of $G$. This space is $X=\mathbb A^1$ for $G=\Gm$ in Riemann--Iwasawa--Tate theory, $X = \GL_2$ for $G= \GL_1\times \GL_2$ in Hecke--Jacquet--Langlands theory, and $X= V \times^{\GL_2^\diag} (\GL_2\times \GL_2)$, where $V$ is the standard representation of $\GL_2$, in Rankin--Selberg theory.
 \item A suitable space of ``Schwartz functions'' $\mathcal F(X(k_v))$ for every completion $k_v$ of $k$; at almost every place, it contains a distinguished vector $\Phi_{0,v}$, giving rise to a restricted tensor product $\mathcal F(X(\mathbb A)) = \bigotimes'_v \mathcal F(X(k_v))$. When $X$ is smooth and $k_v$ is non-Archimedean with ring of integers $\mathfrak o_v$, we have $\Phi_{0,v} = 1_{X(\mathfrak o_v)}$. 
 \item The $X$-theta series
 \[ \Theta: \mathcal F(X(\mathbb A)) \to C^\infty([G]),\]
 where $[G] = G(k)\backslash G(\mathbb A)$, given by $\Theta_\Phi(g):= \Theta(\Phi)(g) := \sum_{\gamma\in X(k)} \Phi(\gamma g)$. This generalizes the Jacobi theta series used by Riemann, and many other series of classical analytic number theory, such as Poincar\'e series (if we allow $X$ to stand for the Whittaker model, which is not just a space but also a non-trivial ``line bundle'' over it, see \S\ref{excellent}), and Eisenstein series (after we pair a suitable theta series with an automorphic form for some Levi subgroup).
\end{itemize}

 The theta series (for varying inputs $\Phi$) are then integrated against automorphic forms $f$, and, under some assumptions on the space $X$, the ``period pairing'' 
 \begin{equation}\label{periodpairing} \langle f, \Theta_\Phi \rangle := \int_{[G]} f(g) \Theta_\Phi(g) \, dg
 \end{equation}
 is expected to be related to a special value of an $L$-function of $f$. This relation will be discussed in \S\ref{conj-global}--\ref{Lconjecture}.

\subsubsection{} \label{singular-intro}

While it is not the main focus of the present article, it should be mentioned that the main point of the proposal of \cite{SaRS} was to include \emph{singular} affine spherical varieties, in which case the ``basic function'' of $X(k_v)$ is the ``IC function,'' obtained through the sheaf--function dictionary from the intersection complex of a suitable geometric model of $X(\mathfrak o_v)$, see \S\ref{singular}. This was inspired by work of Braverman--Kazhdan on the basic affine space \cite{BK1, BK2}, which goes back to the geometric Langlands program \cite{BFGM} and ideas of Drinfeld. 

The conjecture was refined by Ng\^o \cite{Ngo-PS} for a class of affine embeddings of reductive groups; the IC function, for non-Archimedean local fields in equal characteristic, was defined in \cite{BNS}, where Ng\^o's conjecture was proven. In recent joint work with Jonathan Wang \cite{SaWang}, we have obtained similar results for the IC function of a broad class of spherical varieties, including a straightforward generalization of the Hecke and Rankin--Selberg integral represesentations to the Langlands $L$-function associated to the $n$-fold tensor product representation of the dual group
\[ \check G=\underset{n\mbox{ times}}{\underbrace{\GL_2 \times \cdots \times \GL_2}} \xrightarrow\otimes \GL_{2^n}.\]

\subsubsection{}

The elephant in the room, of course, is the global functional equation (and meromorphic continuation), which is not available for these $L$-functions yet. In favorable cases, it should arise from a Poisson summation formula for a ``Fourier transform''
\[ \mathcal F(X(\mathbb A))\to \mathcal F(X^*(\mathbb A)),\]
where $X^*$ is the same variety $X$, with the $G$-action twisted by a Chevalley involution. Such a Fourier transform and a Poisson summation formula are often available for smooth affine spherical varieties, which are vector bundles over homogeneous spaces, but are quite mysterious in the singular case. For the moment, they are known for spaces of the form $X=$ the affine closure of $[P,P]\backslash G$, where $P\subset G$ is a parabolic, by the work of Braverman--Kazhdan \cite{BK1,BK2} (and its refinement \cite{GL2}). An extension to $X=$ the affine closure of $U_P\backslash G$, where $U_P$ is the unipotent radical of $P$, would give rise to the functional equation of normalized Eisenstein series, greatly simplifying and generalizing the theory of $L$-functions obtained through the Langlands--Shahidi method \cite{Langlands-Euler,Shahidi}. In recent work, Getz and his collaborators \cite{GL1,GHL} have proven a Poisson summation formula for a singular space $Y$ which is not directly related to Eisenstein series -- the only example of this sort to date, in my knowledge.

In general, this ``Fourier'' transform may only be available at the level of trace formulas -- see \cite{SaBE2, Xue}, as well as the discussion of \S\ref{Hankel} below.

\subsubsection{} \label{coisotropic}
The ``period pairing'' \eqref{periodpairing} between theta series coming from spherical varieties and automorphic forms is not general enough to include all known integral representations of $L$-functions. At the very least, we need to replace the Schwartz space of a spherical $G$-variety by more general \emph{quantizations of Hamiltonian $G$-spaces}. In the smooth case, those are affine symplectic $G$-spaces $M$, equipped with a moment map $M\to \mathfrak g^*$, which generalize the cotangent bundle $T^*X$ of a smooth spherical $G$-variety. The analog of the ``spherical'' condition for a Hamiltonian $G$-space $M$ is that it be \emph{coisotropic}: namely, that the Poisson algebra $k(M)^G$ of $G$-invariant rational functions on $M$ be Poisson-commutative. 

An example of such a space, that is not the cotangent bundle of a spherical variety, is a symplectic vector space $M$ under the action of a Howe dual pair $G$; that is, $G$ is, up to central isogeny, equal to a product $G_1 \times G_2$ of two subgroups of $\Sp(M)$, where $G_1$ is the commutator of $G_2$, and vice versa. As ``quantization'' of $M$ we understand the Weil representation of the metaplectic group $\operatorname{Mp}(M)$ associated with an additive character $\psi$, restricted to (the metaplectic cover of) $G$. Theta series and the pairing \eqref{periodpairing} still make sense in this setting. 
More general examples mixing the Weil representation with periods are contained in the influential conjectures of Gan, Gross, and Prasad \cite{GGP}. 

In ongoing work with Ben-Zvi and Venkatesh, we describe a class of coisotropic Hamiltonian spaces M whose ``quantizations'' in the form of theta series are expected to be related to special values of $L$-functions, and we demonstrate, by means of known examples, that the $L$-value associated to such a space gives rise to a \emph{dual Hamiltonian space $\check M$ for the Langlands dual group}. In the context of the geometric Langlands program, this leads to a hierarchy of conjectures, with connections to mathematical physics.
We will encounter one of these conjectures in our discussion of unramified $L$-factors in \S\ref{derivedSatake} below.

\subsection{Outline of this paper}

In Section \ref{sec:relative} we introduce the relative Langlands program, up to the conjectural Euler factorization of the period pairings \eqref{periodpairing}.

In Section \ref{sec:Lfunction} we discuss the relationship between the local unramified Euler factors and special values of $L$-functions. 

In Section \ref{sec:BE} we discuss the ``beyond endoscopy'' approach to functoriality, generalized to the setting of the relative Langlands program.

Section \ref{sec:transfer} provides a new interpretation for the transfer operators of functoriality studied in \cite{SaRankone}, based on the concept of quantization. \withoutproofs{(Proofs for the results of this section will appear in an expanded version of this article on the arXiv.) }

Finally, in Section \ref{sec:future} we discuss interesting research directions for the near future.

\subsection{Notation and language} 

\begin{itemize}[leftmargin=1em]
 \item In general, when a variety is defined over a local field $F$, and there is no danger of confusion, we will use the same letter to denote its $F$-points, e.g., ``a Schwartz function on $X$'' really means ``on $X(F)$.''
 \item For a quasiaffine $G$-variety over a field $F$, we will denote by $X/G$ the stack quotient, and by $X\sslash G$ the invariant-theoretic quotient $\operatorname{Spec}F[X]^G$. 
 \item A ``complex line bundle'' on the points of a smooth variety $X$ over a local field $F$ will be
\begin{itemize}[leftmargin=2em]
 \item when $F=\mathbb R$ or $\mathbb C$, a complex line bundle on $X(F)$, viewed as a smooth (Nash) manifold;
 \item when $F$ is non-Archimedean, a locally constant sheaf of complex vector spaces ($l$-sheaf) on $X(F)$, for the $p$-adic (Hausdorff) topology, with $1$-dimensional stalks.
\end{itemize}
When no confusion arises, we will just say ``line bundle'' for a complex line bundle; when we want to distinguish it from a line bundle on $X$ in the sense of algebraic geometry, we will say ``algebraic line bundle'' for the latter. 
\item An algebraic line bundle $L$ over a smooth $F$-variety $X$, where $F$ is a local field, together with a complex number $s$, give rise to a complex line bundle $|L|^s$ on $X(F)$, by reduction of the corresponding $\Gm$-torsor via the sequence of maps $\Gm(F)\xrightarrow{|\bullet|} \mathbb R^\times_+ \xrightarrow{ x\mapsto x^s} \mathbb C^\times$. (When $F$ is non-Archimedean, the absolute value map is discretely-valued,  giving rise to the structure of a locally constant sheaf.)
 \item When $L = \det T^*X$, the line bundle of volume forms on $X$, the associated complex vector bundle $|L|^s$  is known as the bundle of $s$-densities. We will, in general, understand the field $F$ as endowed with a Haar measure; this identifies densities, i.e., sections of $|L|$, as measures on $X(F)$. When no confusion arises, we will denote $|dx|$, the density attached to a volume form $dx$, simply by $dx$.
 \item The space of Schwartz functions on the $F$-points of a variety $X$, where $X$ is a local field, will be denoted by $\mathcal F(X(F))$, the space of Schwartz measures by $\mathcal S(X(F))$, and the space of Schwartz half-densities by $\mathcal D(X(F))$. These are smooth, compactly supported sections of the corresponding bundles of $s$-densities, in the non-Archimedean case. In the Archimedean case, they are smooth sections of rapid decay, see \cite{AGSchwartz}. We will also say ``test functions/measures,'' etc., for ``Schwartz.''
 
 \item For an admissible, smooth, complex representation $\pi$ of a reductive group over a local field, we will denote by $\tilde\pi$ its contragredient. When $\pi$ is unitary, $\tilde\pi$ is identified with the complex conjugate $\bar\pi$.
 \item We will generally prefer to replace a hermitian pairing $H$ between functions by the associated bilinear pairing $B(\Phi_1, \Phi_2)=H(\Phi_1,\overline{\Phi_2})$. When $H$ is an inner product, we will sometimes call $B$, by abuse of language, an ``inner product.'' 
 
 \item $\mathcal W_F$ will denote the Weil group of a local or global field, and $\mathcal L_F$ will be the ``Langlands group,''  whose representations should  parametrize local and automorphic $L$-packets.  It is the Weil group for Archimedean local fields and global function fields, the Weil--Deligne group for non-Archimedean local fields, and a conjectural extension of the Weil group for number fields. 
 
 \item We adopt the ``Weil group'' convention for $L$-groups of reductive groups, ${^LG} = \check G\rtimes \mathcal W_F$; the dual group $\check G$ is identified with the set of its complex points.
 
\end{itemize}

\section{The relative Langlands conjectures} \label{sec:relative}

\subsection{The local and global spectrum of a spherical variety}

\subsubsection{} To understand the relationship between period pairings \eqref{periodpairing} and $L$-functions, one needs to start by understanding the phenomenon of \emph{distinction}, highlighted by the groundbreaking work of Jacquet and his collaborators \cite{JLR, Jacquet-symmetric}. A naive formulation of this phenomenon goes as follows:

\begin{quote}
 The local and global spectrum of a spherical $G$-variety $X$ only contain representations with Langlands parameters in a certain subgroup $\LG_X\subset \LG$ of the $L$-group of $G$. 
\end{quote}

To make sense of this statement, we need to explain ``the local and global spectrum of a $G$-variety.'' Then, we need to talk about the $L$-group $\LG_X$. Finally, the statement needs to be corrected, for some ``non-tempered'' varieties $X$, replacing Langlands parameters by appropriate Arthur parameters.

Let $R$ denote either a local field, or the adelic points of a global field. 
In order to define the local and global spectrum of a $G$-variety $X$ (defined over the corresponding field), we will introduce the Plancherel formula and the relative trace formula. These decompose certain distributions -- or rather, generalized functions -- on the $R$-points of $X\times X$, invariant under the diagonal action of $G$. For the purposes of the Langlands program, it turns out to be more natural to think of them as generalized functions on the $R$-points of the quotient stack $\mathfrak X = (X\times X)/G^\diag$, which naturally includes ``pure inner forms'' of the pair $(G,X)$.

\subsubsection{} \label{localspectrum}

Let $F$ be a local field. 
The space $L^2(X)$ is the Hilbert space completion of the space $\mathcal D(X)$ of Schwartz half-densities on $X(F)$, with respect to the $L^2$-inner product, and furnishes a unitary representation of $G$. 
By the Plancherel decomposition, there are a measure $\mu_X$ on the unitary dual $\widehat G$ of $G$, and a measurable family of linear forms
\[J_\pi: \mathcal D(X\times X)\to \mathbb C,\]
\begin{minipage}[t]{0.11\linewidth}
such that:
\end{minipage} 
\begin{minipage}[t]{0.89\linewidth}
\begin{itemize}[leftmargin=2em]
 \item for $\mu_X$-almost every $\pi$, $J_\pi$ factors as $\mathcal D(X\times X) \to \pi\hat\otimes \bar\pi \to \mathbb C$, and
 \item for all $\Phi\in \mathcal D(X\times X)$, we have 
\end{itemize}
\end{minipage}
\begin{equation}\label{Plancherel} \int_{X^\diag} \Phi = \int_{\hat G} J_\pi(\Phi) \mu_X(\pi).
\end{equation}

A linear form satisfying the first property above will be called a \emph{relative character}.
The product $J_\bullet\mu_X$, which can be thought of as a measure valued in the space of functionals on $\mathcal D(X\times X)$, is uniquely defined. Moreover, the relative characters $J_\pi$ are invariant under the diagonal action of $G=G(F)$; thus, they factor through the coinvariant space $\mathcal D(X\times X)_G=$ the quotient of $\mathcal D(X\times X)$ by the (closed, in the Archimedean case) subspace generated by elements of the form $f-g\cdot f$, where $g\cdot f$ denotes the action of $g\in G$ on $f$ by diagonal translation.

Let us assume that $X$ carries a positive $G$-invariant measure $dx$, and use it to identify functions, half-densities, and measures on $X$ through the $G$-equivariant maps $\Phi\mapsto \Phi (dx)^\frac{1}{2} \mapsto \Phi dx$ (and similarly on $X\times X$). Then, the coinvariant space $\mathcal D(X\times X)_G \simeq \mathcal S(X\times X)_G$ is more naturally understood as a subspace of the \emph{Schwartz space of the quotient stack $\mathfrak X:= (X\times X)/G$} \cite{SaStacks}. This Schwartz space is really a complex of vector spaces, but here we will focus only on its zeroth cohomology, which has the explicit description
\begin{equation}\label{Schwartz}
 \mathcal S(\mathfrak X) = \bigoplus_{\alpha} \mathcal S(X^\alpha \times X^\alpha)_{G^\alpha}.
\end{equation}
Here, $\alpha$ runs over isomorphism classes of $G$-torsors (parametrized by the Galois cohomology set $H^1(\Gamma_F, G)$, where $\Gamma_F$ is the Galois group of a separable closure of $F$); if $R^\alpha$ is a representative of a class $\alpha$, we let $G^\alpha  = \operatorname{Aut}_G(R^\alpha)$, and $X^\alpha = X\times^G R^\alpha$, a $G^\alpha$-space. In other words, $G^\alpha$ is what is called a ``pure inner form'' of $G$, and $X^\alpha$ can similarly be called a ``pure inner form'' of $X$, if its set of $F$-points is nonempty. 

The Plancherel formula \eqref{Plancherel} extends to $\mathcal S(\mathfrak X)$, with a measure $\mu_{\mathfrak X}$, on the right hand side, on the union of the unitary duals of the pure inner forms $G^\alpha$. The support $\Pi_{\mathfrak X}$ of this measure (avoiding redundancy -- i.e., the support of the canonical linear form-valued measure $J_\bullet \mu_{\mathfrak X}$) can be called the \emph{local ($L^2$-)spectrum of the quotient stack $\mathfrak X$}.

\subsubsection{} \label{globalspectrum} The global (automorphic) spectrum of $X$ (or rather, again, of the stack $\mathfrak X = (X\times X)/G^\diag$) can be defined through the relative trace formula of Jacquet. This is a generalization of the Arthur--Selberg trace formula, and an automorphic analog of the local Plancherel formula. Its definition uses the theta series encountered in \S\ref{thetaseries}, therefore we assume here that $X$, defined over a global field $k$, is quasiaffine, so that $X(k)$ is discrete in the adelic points $X(\mathbb A)$. As before, we write $[G]=G(k)\backslash G(\mathbb A)$ for the automorphic quotient space.  

Roughly speaking, the relative trace formula is the Plancherel formula for $L^2([G])$, applied to the inner product of two theta series for $X$, i.e., decomposing the functional 
\begin{equation}\label{RTF-naive} \RTF_{X}: \mathcal F(X(\mathbb A))\otimes \mathcal F(X(\mathbb A)) \ni \Phi_1\otimes \Phi_2 \mapsto \Theta_{\Phi_1} \otimes \Theta_{\Phi_2} \mapsto \int_{[G]} \Theta_{\Phi_1}(g) \Theta_{\Phi_2}(g) \, dg \in \mathbb C.
\end{equation}

This naive point of view requires some caution:
\begin{itemize}[leftmargin=2em]
\item The inner product on the right hand side of \eqref{RTF-naive} does not, in general, converge, and needs to be regularized. Depending on $X$, there may be a \emph{canonical} way to regularize it, described in \cite[\S6]{SaStacks}. In many cases of interest, though, notably in the case of the Arthur--Selberg trace formula (where $X=H$, a reductive group, and $G=H\times H$), a canonical regularization is not available, and it takes the mastery of Arthur's work \cite{Arthur-invariant} to engineer an invariant expression. Such work has not yet been done in the general setting of the relative trace formula.
 
 \item We can again choose a $G$-invariant measure on $X(\mathbb A)$ (e.g., Tamagawa measure) to identify functions with measures, and understand the $G(\mathbb A)^\diag$-invariant functional \eqref{RTF-naive} as a functional on 
  $\mathcal S(X\times X(\mathbb A))_{G(\mathbb A)}.$
 As in the local case, this space is a subspace of the global Schwartz space of the stack $\mathfrak X= (X\times X)/G^\diag$, 
 \[ \mathcal S(\mathfrak X(\mathbb A)) = \bigotimes'_v \mathcal S(\mathfrak X(k_v)),\]
 and the relative trace formula should be defined as a functional on the bigger space,
\[ \RTF_{\mathfrak X} = \sum_\alpha \RTF_{X^\alpha},\]
where now $\alpha$ runs over isomorphism classes of $G$-torsors over the global field $k$. 
\end{itemize}

Ignoring the regularization issue, if we could apply the Plancherel formula for $\bigoplus_\alpha L^2([G^\alpha])$ to the pairing \eqref{RTF-naive}, we would obtain the \emph{spectral side} of the relative trace formula,
\begin{equation}\label{RTF-spectral} \RTF_{\mathfrak X} = \int J^\aut_{\pi} \mu_{\mathfrak X}^\aut(\pi), 
\end{equation}
where the product $J^\aut_\bullet \mu_{\mathfrak X}^\aut$ is a measure on the union $\bigsqcup_\alpha \widehat{G^\alpha}^\aut$ of ($L^2$-)automorphic spectra of the pure inner forms of $G$, valued in linear forms on $\mathcal S(\mathfrak X(\mathbb A))$. 

The global (automorphic) spectrum of $\mathfrak X$ is defined as the support of $J^\aut_\bullet \mu_{\mathfrak X}^\aut$. Clearly, this definition is incomplete, as it relies on overcoming the aforementioned issues of regularization, and developing a spectral decomposition for the relative trace formula.

\subsection{The Langlands dual group}

\subsubsection{} The local and global spectrum of a spherical variety $X$ are conjecturally governed by the \emph{$L$-group} $\LG_X$ of $X$. We owe this dual group to the insights developed by Nadler in his thesis \cite{Na-real}, and in his joint work with Gaitsgory \cite{GN}. They realized that the ``little Weyl group'' of a spherical $G$-variety (defined by Brion in \cite{Brion}, and generalizing the little Weyl group of a symmetric space) corresponds to a subgroup $\check G_X\subset \check G$ of the Langlands dual group of $G$, and gives rise to a form of the geometric Satake isomorphism for the spherical variety. In \cite{SV}, it was proposed that this dual group comes equipped with a distinguished morphism
\begin{equation}\label{dualgroup}
 \check G_X\times \SL_2 \to \check G
\end{equation}
that governs the harmonic analysis of $X$, in a way that will be described below. Since Gaitsgory and Nadler did not fully identify their dual group $\check G_X$ (constructed in a Tannakian way), an independent description of a morphism \eqref{dualgroup} was achieved by Knop and Schalke \cite{KnSch}; we can take this as the definition of the dual group, for what follows. Finally, for the purposes of the Langlands program, we need an $L$-group, in the form of an extension 
\[ 1\to \check G_X \to {^LG_X} \to \mathcal W_F \to 1.\]
The correct definition of this $L$-group, when $G$ is not split, is not completely understood yet, although it is probably within reach. 
For what follows, we will assume such an $L$-group, and an extension of the homomorphism \eqref{dualgroup} to the $L$-groups, in the sense that the conjectures to be stated should hold for an appropriate definition of ${^LG_X}$. 

\subsubsection{} \label{invariants} We briefly describe one way to characterize the root datum of the dual group $\check G_X$: As in the case of reductive groups, the first step is to describe a canonical maximal torus, which in turn is dual to an ``abstract Cartan'' group. Let $A$ be the abstract Cartan group of $G$; it is canonically equal to the reductive quotient of any Borel subgroup of $G$. Fix such a Borel subgroup $B$, with unipotent radical $N$, and let $X^\circ$ be the open $B$-orbit. On the quotient $X^\circ \sslash N$, $B$ acts through a quotient $A_X$ of $A$; this is \emph{the Cartan group of $X$}, and it can be seen to be independent of $B$, in the sense that any two choices induce canonical tori up to a canonical isomorphism. (These definitions assume that $B$ is defined over the base field, but by Galois descent the Cartan groups $A$ and $A_X$ are defined over the field, even if $B$ is not.) 

The quotient $A\to A_X$ gives rise to a morphism of dual tori $\check A_X\to \check A$, which could have nontrivial (finite) kernel. The image of this morphism is the canonical maximal torus of the Gaitsgory--Nadler dual group $\check G_X$. (We caution the reader that in \cite{SV} the group $\check G_X$ was not necessarily defined as a subgroup of $\check G$, and had $\check A_X$ as its maximal torus.)

\label{sphericalroots}
It is slightly harder to define the little Weyl group $W_X$. Once this is done, the coroots of $\check G_X$, which will be called the \emph{normalized spherical roots} of $X$, are uniquely determined up to multiple, and that multiple is fixed by the following axiom:
\begin{quote}
 A normalized spherical root is either a root of $G$, or the sum of two strongly orthogonal roots, i.e., two roots whose linear span contains no other roots but their multiples.
\end{quote}

\subsubsection{} \label{Knopwork} For the Weyl group, there are many equivalent definitions. Most relevant to our purposes, when $X$ is defined over a field $F$ in characteristic zero, is the following one, due to Knop \cite{KnWeyl}: We may assume that $X$ is smooth and that $F$ is algebraically closed (since $\check G_X$ only depends on the open $G$-orbit over the algebraic closure). Consider the cotangent space $M=T^*X$, equipped with the moment map $\mu: M\to \mathfrak g^*$. If $\mathfrak a^*$ denotes the dual Lie algebra of the Cartan of $G$, Chevalley's isomorphism identifies the invariant-theoretic quotient $\mathfrak g^*\sslash G$ with $\mathfrak a^*\sslash W$. The \emph{polarized cotangent bundle}
\[ \widehat{M} := M\times_{\mathfrak a^*\sslash W} \mathfrak a^*\]
is not, in general, irreducible. Knop describes a distinguished irreducible component $\widehat{M}^\circ$ living over the dual Lie algebra $\mathfrak a_X^*\subset \mathfrak a^*$ of $A_X$, and shows that the map $ \widehat{M}^\circ \to M$ is generically a Galois cover with covering group a subquotient $W_X$ of the Weyl group; this is the little Weyl group of $X$ \cite[\S6]{KnWeyl}.

For later use, we mention a related result of Knop, still in the homogeneous case: Let $\mathfrak g_X^*=$ the normalization of the image of the moment map \emph{in $M$} (i.e., the spectrum of the integral closure of the image of $F[\mathfrak g^*]$ in $F[M]$). The composition $\widehat{M}^\circ \to \mathfrak a_X^*\to \mathfrak c_X^*:= \mathfrak a_X^*\sslash W_X$ factors through a map $\mu_G: M \to \mathfrak g_X^*\to \mathfrak c_X^*$, called the \emph{invariant moment map}, and identifies $\mathfrak c_X^*$ with the invariant-theoretic quotient $M\sslash G$ \cite[Korollar 7.2]{KnWeyl}. 

\subsubsection{} \label{Arthur} Finally, the restriction of the map \eqref{dualgroup} to $\SL_2$, which we will call the ``Arthur-$\SL_2$'' of $X$, is determined by the conjugacy class of parabolics of the form
\[ P(X) = \{ g\in G| X^\circ g = X^\circ\},\]
where $X^\circ$ is the open orbit for a Borel subgroup $B$. In the quasiaffine case, $P(X)$ is the largest parabolic such that all highest weight vectors in $k[X]$ are $P(X)$-eigenvectors. To this class is canonically associated a standard Levi subgroup $\check L(X)$ of $\check G$, and the Arthur-$\SL_2$ of $X$ is a principal $\SL_2 \to \check L(X)$.

\subsubsection{} \label{excellent} In order to keep the discussion that follows as simple as possible, let us single out a convenient class of spherical varieties: We will say that a spherical $G$-variety $X$ is \emph{excellent} if it is affine, homogeneous, and the kernel of the map $A\to A_X$ is connected (equivalently, the map $\check A_X\to \check A$ of dual tori is injective). 

We also need to enlarge the class of spherical varieties, in order to include objects such as the Whittaker model. The Whittaker model is the space $N\backslash G$, where $G$ is quasisplit and $N$ is a maximal unipotent subgroup, endowed with a nondegenerate character $\psi: N(F)\to \mathbb C^\times$. This character defines, by induction, a complex line bundle $L_\psi$ over $N\backslash G$, and the Whittaker model consists of sections of this line bundle. In the sequel, when we say that $Y$ is ``the Whittaker model,'' we will mean the space $N\backslash G$ together with this line bundle, and we will be using the Schwartz space notation $\mathcal F(Y)$, $\mathcal S(Y)$, etc., to denote Schwartz sections (resp.\ measures) valued in this line bundle. For a more general discussion of ``Whittaker induction,'' see \cite[\S2.6]{SV}.

\subsection{Conjectures}

\subsubsection{} \label{conj-local} Let $X$ be defined over a local field $F$, and let $\Pi_{\mathfrak X}$ be the set of $L^2$-distinguished representations of $X$ and its pure inner forms, as in \S\ref{localspectrum}. We assume, for simplicity, that $X$ carries an invariant measure, to identify measures with half-densities.  The ``relative local Langlands conjecture'' of my work with Venkatesh \cite[\S16]{SV} states:
\begin{conjecture*}
Let $\mu_{{^LG_X}}$ be the  natural measure on the set of tempered local Langlands parameters into ${^LG_X}$.
There is a decomposition of the inner product on $\mathcal D(X\times X)$,
\begin{equation}\label{Plancherel-packets} \int_{X^\diag} \Phi = \int J_\phi(\Phi) \mu_{{^LG_X}}(\phi),
\end{equation}
where the ``stable relative characters'' $J_\phi$ are linear combinations of relative characters for representations belonging to Arthur packets with parameter  
\begin{equation}\label{XArthur}\mathcal L_F \times \SL_2 \xrightarrow{\phi \times \operatorname{Id}} {^LG_X}\times\SL_2 \xrightarrow{\eqref{dualgroup}} {^LG}.
\end{equation}
\end{conjecture*}
For the ``natural measure'' on such parameters, see \cite[\S16]{SV}. Comparing with the Plancherel formula \eqref{Plancherel}, the conjecture implies that the local $L^2$-spectrum $\Pi_{\mathfrak X}$ of $\mathfrak X$ belongs to the union of Arthur packets with parameters of the form \eqref{XArthur}. Developing a Plancherel formula for $X$ in terms of discrete-mod-center spectra of its ``boundary degenerations,'' as in \cite{HC3, Waldspurger-Plancherel, Delorme-real, vdBS1, vdBS2, SV, Delorme-Plancherel}, reduces the conjecture to discrete spectra.

When $G$ is quasisplit, acts faithfully on $X$, and the map $\check A_X\to \check A$ is injective (\S \ref{invariants}), one would expect the functionals $J_\phi$ of \eqref{Plancherel-packets}, after summing over all pure inner forms of $X$, to be nonzero. In a broad range of individual cases, including the Gan--Gross--Prasad conjectures \cite{GGP} and other cases considered by D.\ Prasad \cite{Prasad-relative} and C.\ Wan \cite{Wan}, we have much more precise conjectures about how many and which elements in the given Arthur packets are distinguished. A number of cases have been proven by Waldspurger, M{\oe}glin, Beuzart-Plessis, Gan, Ichino and others \cite{Waldspurger-GP, MW-GP, Beuzart-GP, Gan-Ichino}, and by M{\oe}glin--Renard \cite{MR} for symmetric spaces over $\mathbb R$.

Besides the question of $L^2$-distinction, one can ask the question of smooth distinction: which irreducible representations embed as $\pi\hookrightarrow C^\infty(X)$? The general answer to this question is less understood.

\subsubsection{} \label{conj-global} Now, let $X$ be defined over a global field $k$, and let $\Pi_{\mathfrak X}^\aut$ be the automorphic spectrum of $\mathfrak X$, as in \S\ref{globalspectrum}. We recall that its definition is, in general, conditional on developing the spectral decomposition of the relative trace formula. Nonetheless, one can often restrict to parts of the spectrum where the full relative trace formula is not needed; for example, if $\pi$ is a discrete automorphic representation where the period pairing \eqref{periodpairing} is absolutely convergent, the corresponding functional-valued measure $J^\aut_\bullet \mu^\aut_{\mathfrak X}$ of \eqref{RTF-spectral}, applied to a test function $\Phi_1\otimes \Phi_2$ on $(X\times X)(\mathbb A)$, should have the meaning of 
\[J^\aut_\pi (\Phi_1\otimes\Phi_2) \mu^\aut_{\mathfrak X}(\pi) =  \sum_{f} \Big(\int_{[G]} \Theta_{\Phi_1} f\Big) \Big(\int_{[G]} \Theta_{\Phi_2} \bar f\Big),\]
where $f$ runs over an orthonormal basis of $\pi$.

A landmark in our understanding of these global relative characters was the paper \cite{II} of Ichino and Ikeda, generalizing the formula of Waldspurger \cite{Waldspurger} to a precise conjectural Euler factorization of $J^\aut_\pi$, in the case of orthogonal Gross--Prasad periods. The conjectures of Ichino and Ikeda gave rise to the realization that there was a general pattern in the Euler factorization of automorphic periods, and were quickly adapted to other cases. Unlike the orthogonal case, which remains open, the conjecture for unitary Gross--Prasad periods has been proven in \cite{Zhang,BLZZ,BCZ}, its analog for Whittaker periods of metaplectic and unitary groups was proven in \cite{LM, LM-unitary, Morimoto}, and there are significant partial results in many other cases. 

A generalization of the Ichino--Ikeda conjecture to a wide range of spherical periods (satisfying certain conditions) was proposed in \cite{SV}. As in the local case, it lacks the precision of conjectures known in special cases, hence leaving an open problem that should be addressed in the near future. On the other hand, the conjecture of \cite{SV} makes clear the connection between the (global) relative trace formula and the (local) Plancherel formula. I will formulate a variant of this conjecture here, using the hypothetical notion of global Arthur parameters, and being a bit vague on choices of measures (see \cite[\S17]{SV} for some hints). Its formulation also relies on Conjecture \ref{Lconjecture} below, expressing the local Plancherel density of the basic function $\Phi_{0,v} \in \mathcal F(X(k_v))$ at almost every place $v$ in terms of a local $L$-value $L_X(\phi_v):= L(\phi_v, \rho_X, 0)$, where $\phi_v$ is a local unramified Langlands parameter into ${^LG_X}$ and $\rho_X: {^LG_X}\to \GL(V_X)$ is a certain representation of the $L$-group of $X$.

\begin{conjecture*}
 There is a decomposition 
 \[\RTF_{\mathfrak X}  =  \int J_\phi^\aut \mu_X^\aut(\phi),\]
 where $\mu_X^\aut$ is a measure on the set of global Arthur parameters which factor as
 \[ \mathcal L_k \times \SL_2 \xrightarrow\phi {^LG_X}\times \SL_2 \xrightarrow{\eqref{dualgroup}} {^LG}\]
 (with $\phi$ lying over the identity map for the projections to $\SL_2$), and $J^\aut_\phi$ is a sum of relative characters $\mathcal S(\mathfrak X(\mathbb A))\to (\pi\hat\otimes\bar\pi)_{G(\mathbb A)}\to \mathbb C$ for automorphic representations $\pi$  belonging to the corresponding Arthur packet. 
  
 Moreover, when $\mathfrak X$ is stable, the restriction of $J_\phi^\aut \mu_X^\aut(\phi)$ to the \emph{most tempered} Arthur type ($\phi(\SL_2) = \SL_2$), \emph{away from the poles of $L_X(\phi|_{\mathcal L_k})$} is equal to 
 \begin{equation}\label{Euler} \frac{1}{|S_\phi|} \prod_v' J_{\phi_v} \cdot \mu_{{^LG_X}}(\phi),
 \end{equation}
 where $\mu_{{^LG_X}}$ is the natural measure on the set of such parameters, $S_\phi$ is the stabilizer of $\phi$ in $\check G_X$, and
 the factors $J_{\phi_v}$ of the Euler product are the local Plancherel relative characters of Conjecture \ref{conj-local}.
\end{conjecture*}

``Stable,'' here, means that the stabilizers of generic points have trivial Galois cohomology; one can drop this assumption, replacing $\RTF_{\mathfrak X}$ by its (properly defined) stable analog. The Euler product of the conjecture needs to be understood, outside of a finite set $S$ of places, as the partial $L$-value ${L_X^S(\phi)}/{L^S(\phi, \check{\mathfrak g}_X, 1)}$, according to Conjecture \ref{Lconjecture} below. The conjecture can be generalized to other quantizations of suitable Hamiltonian spaces, such as the theta correspondence, where it was shown in \cite{SaHowe} to follow from a version of the Rallis inner product formula proven in \cite{GQT, Yamana}. The conjecture is compatible with earlier results and methods for computing period integrals, such as the ``unfolding'' method \cite[\S18]{SV}, or the work of Jacquet and Feigon--Lapid--Offen on unitary periods \cite{Jacquet-factorization,FLO,Beuzart-unitary}.

\section{The $L$-value of a spherical variety} \label{sec:Lfunction}

\subsection{Plancherel density of the basic function}

\subsubsection{} \label{Lconjecture}
It is a very interesting problem to relate the Euler factors of \eqref{Euler}  -- that is, the local Plancherel densities -- to special values of local $L$-functions at every place, including ramified and Archimedean ones. However, we will confine ourselves here to the calculation of the local Plancherel density of the \emph{basic function} $\Phi_0 \in \mathcal F(X(F))$, for a local non-Archimedean field $F$. We assume that $G, X$ are defined over the integers $\mathfrak o$ of $F$, with $G$ reductive, and recall that the basic function is equal to $1_{X(\mathfrak o)}$, when $X$ is smooth and affine; in general, it is the ``IC function,'' see \S\ref{singular} below. We assume that the map $\check A_X\to \check A$ is injective (\S\ref{invariants}). 

For simplicity of presentation, we will assume that $G$ is split, so that the maximal compact subgroup $\check A_X^1\subset \check A_X$ is identified with the group of unramified unitary characters of $A_X$. The unramified representations appearing in Conjecture \ref{conj-local} are those obtained by unitary induction of those characters from the parabolic $P(X)$ (\S\ref{Arthur}) through the quotient $P(X)\to A_X$, and the ``natural measure'' of the conjecture, restricted to unramified parameters (i.e., to $\check A_X^1/W_X$), reads
\[ \mu_{\check G_X}(\phi) = \frac{L(\phi, \check{\mathfrak g}_X/\check{\mathfrak a}_X, 1)}{L(\phi, \check{\mathfrak g}_X/\check{\mathfrak a}_X, 0)} d_{\operatorname{Haar}}\phi.\]

\begin{conjecture*}
 Let $X$ be an affine spherical variety satisfying the conditions above, with a good model over $\mathfrak o$, and $\Phi_0$ its basic function. There is a representation $\rho_X:{^LG_X}\to \GL(V_X)$ such that, setting $L_X(\phi)=L(\phi, \rho_X, 0)$, the Plancherel decomposition of $\Phi_0$ reads
 \begin{equation}\label{density-conj}
  \Vert \Phi_0\Vert^2 = \int_{\check A_X^1/W_X} \frac{L_X(\phi)}{L(\phi, \check{\mathfrak g}_X, 1)} \, \mu_{\check G_X}(\phi).
 \end{equation}
\end{conjecture*}

We refrain from giving details on the precise normalization of Haar measures, or the precise meaning to ``good model;'' at a minimum, the conjecture should be valid at almost every place for any global model. This, in particular, will identify almost every Euler factor of Conjecture \ref{conj-global} as a quotient of special values of local $L$-functions. Note that the ``true'' point of evaluation of $L_X$ is not $0$, but is encoded in $\rho_X$, which is a representation of the full $L$-group. Here, this representation would factor through the unramified quotient, and the ``true'' point of evaluation depends on the action of Frobenius.

The relation between local $L$-values and Plancherel densities is a fascinating one. On the surface, it is just the outcome of a local integral. For example, when $X(\mathfrak o)=H(\mathfrak o)\backslash G(\mathfrak o)$, the value $J_\phi(1_{X(\mathfrak o)})$ is given by the following Ichino--Ikeda local period, in the so-called \emph{strongly tempered} cases where it is convergent:
\[ J_\phi (1_{X(\mathfrak o)}) = \int_H m_\phi(h) dh,\]
where $m_\phi$ is the zonal spherical function (= unramified matrix coefficient with value $1$ at the identity) for the unramified representation with Satake parameter $\phi$.
This calculation, however, has a conceptual meaning, in terms of both harmonic analysis and geometry. We will only attempt to give a flavor of the richness of the topic here.

\subsubsection{} \label{density-smooth}
The study of the Plancherel density of the basic function is a topic with a long history. The mainstream method for calculating it is the Casselman--Shalika method \cite{Casselman,CS}, and it is essentially equivalent to the problem of calculating eigenvectors for the unramified (spherical) Hecke algebra $\mathcal H(G(F), G(\mathfrak o))$ on the space $C^\infty(X(F))$. 

The calculation was related to the structure theory of spherical varieties in \cite{SaSph}, for $G$ split. Here, we will formulate the result in the special case when $X$ is an excellent spherical variety (\S\ref{excellent}) with $\check G_X=\check G$. Fix a Borel subgroup $B\subset G$, with unipotent radical $N$. The important geometric invariants determining the $L$-value are the \emph{colors} of the spherical variety $X$: those are the $B$-stable prime divisors on $X$, over the algebraic closure. For simplicity, we will assume all those divisors to be defined over $F$. 
Each such divisor $D$ induces a valuation on the function field $F(X)$, which we restrict to the multiplicative group of nonzero $B$-eigenfunctions. This gives rise to a homomorphism factoring as 
\[ F(X)^{(B)} \to X^\bullet(A_X) \to \mathbb Z,\]
i.e., an element $\check v_D \in X_\bullet(A_X)$, the character group of the dual torus $\check A_X$ (which here is equal to $\check A$). Assume -- as in all cases I have checked -- that the weights $\check v_D$ are all minuscule.\footnote{The published version of the article cites corollary from \cite{SaWang}, which asserts the minuscule property; however, the argument of that corollary has a gap, and proves a weaker statement.} Let $V_X$ be the smallest self-dual (algebraic) representation of $\check G$ which contains all those weights (with multiplicity, if some of the $\check v_D$'s coincide). For an alternative interpretation of this representation, in terms of the structure of the Hamiltonian space $T^*X$, see Theorem \ref{Hamiltonian} below. The following was proven in \cite{SaSph} under some assumptions, and in \cite{SaWang} in general:

\begin{theorem*}
 The Plancherel density of the basic function $1_{X(\mathfrak o)}$ is given by 
 \begin{equation}
J_\phi(1_{X(\mathfrak o)} \otimes 1_{X(\mathfrak o)}) \mu_{\check G}(\phi)= \frac{L(\phi, V_X, \frac{1}{2})}{L(\phi, \check{\mathfrak g}, 1)}\mu_{\check G}(\phi).
 \end{equation}
\end{theorem*}

Note that, for simplicity, we have assumed that $\check G_X=\check G$. The point $\frac{1}{2}$ of evaluation changes in the general case.

\subsubsection{} \label{singular}
The case of \emph{singular} affine varieties $X$ was undertaken in \cite{SaWang}, for the cases with $\check G_X= \check G$. As mentioned, here one needs to work in a geometric setting, assuming that $F$ is a local field in equal characteristic, $F\simeq  \mathbb F_q((t))$, with $G, X$ defined over $\mathbb F_q$. (There are ad hoc ways to transfer the results to mixed characteristic, but it would be nice to see a direct geometric approach.) The basic function $\Phi_0$ is then defined as the Frobenius trace on the stalks of the intersection complex of finite-dimensional formal models of $L^+X$, the formal arc space of $X$ \cite{BNS}.

Let us discuss the special case when $X$ is the affine closure $\operatorname{Spec} \mathbb F_q[X^\bullet]$ of its open $G$-orbit $X^\bullet$. Colors, here, do not need to be minuscule, but one can still define $V_X$ as before.
We have the following generalization of Theorem \ref{density-smooth}:

\begin{theorem*}
 There is a representation $V_X'$ of $\check A$, with the same weights as $V_X$ and $W$-invariant multiplicities, such that the Plancherel density of the basic function $\Phi_0$ is
 \begin{equation}
J_\phi(\Phi_0\otimes \Phi_0)\mu_{\check G}(\phi)= \frac{L(\phi, V_X', \frac{1}{2})}{L(\phi, \check{\mathfrak g}, 1)}\mu_{\check G}(\phi).
 \end{equation}
\end{theorem*}

Of course, we expect that $V_X'=V_X$. This is automatic in the minuscule case. For example \cite[Example 1.1.3]{SaWang}, there is a family of varieties $X_n$, $n\in \mathbb N$, which gives rise to the generalization of the Hecke and Rankin--Selberg integrals, mentioned in \S\ref{singular-intro}.

\subsection{Derived Satake equivalence for spherical varieties} \label{derivedSatake}

Ongoing joint work with Ben-Zvi and Venkatesh has revealed deeper relations between periods and $L$-functions; currently, those can be formulated over function fields and their completions. In the local setting, $F\simeq\mathbb F_q((t))$, Conjecture \ref{Lconjecture} should be obtained by applying the sheaf-function dictionary to a categorical statement, along the following lines:

We retain the assumptions of the previous subsection, with $X$ and $G$ defined over $\mathbb F_q$, and also assume $X$ to be smooth. We denote formal loop and arc spaces by $L$, resp.\ $L^+$, so that $LX(\mathbb F_q)=X(F)$, $L^+G(\mathbb F_q) = G(\mathfrak o)$.  
For appropriate measures, the left hand side of \eqref{density-conj}, and, more generally, the pairing of two $G(\mathfrak o)$-invariant functions $f$, $g$ obtained 
 via the sheaf-function dictionary from objects $\mathscr F, \mathscr G$ in the bounded derived category $\operatorname{Shv}(LX/L^+G)$ of constructible $\ell$-adic \'etale sheaves on $LX/L^+G$ can be computed as the (alternating) trace of geometric Frobenius on a derived homomorphism complex:
\[ \int_{X(F)} f(x) g(x) dx = \operatorname{tr}\Big( \operatorname{Frob}_q, \operatorname{\emph{Hom}} (\mathscr F, D\mathscr G)^*\Big),\]
where $D=$Verdier dual. The pairing is really a finite sum, and makes sense over $\overline{\mathbb Q_\ell}$.

The right hand side of \eqref{density-conj}, through a simple application of the Weyl integration formula, can be interpreted as the Frobenius trace on 
\[\mathbb C[V_X]^{\check G_X} = \mathbb C[\check M]^{\check G},\]
where $V_X$ is the space of the representation $\rho_X$, and we have set 
\[\check M = V_X\times^{\check G_X} \check G.\]
The empirical observation is that the space $\check M$ \emph{has a natural symplectic structure}, and, moreover, that \emph{the assignment $M=T^*X \to \check M$ is involutive}, although, to make sense of this, one needs to allow for more general coisotropic Hamiltonian spaces, as mentioned in \S\ref{coisotropic}. For the categorical analog of Conjecture \ref{Lconjecture}, we need to \emph{shear} the ring $\mathbb C[\check M]$ into a dg-algebra $\mathbb C[\check M]^\shear$, with zero differentials, in  degrees related to the action of Frobenius in $\rho_X$. 

\begin{conjecture*}
 Fix an isomorphism $\mathbb C = \overline{\mathbb Q_\ell}$. There is an equivalence of triangulated $\mathbb C$-linear categories
 \[ \operatorname{Shv}(LX/L^+G) \xrightarrow\sim D_{\operatorname{per}}^\shear(\check M/\check G),\]
where $D_{\operatorname{per}}^\shear(\check M/\check G)$ denotes the full triangulated subcategory, generated by perfect complexes, of the category of $\check G$-equivariant differential graded $\mathbb C[\check M]^\shear$-modules localized by quasi-isomorphisms.
\end{conjecture*}

This generalizes the derived Satake equivalence of Bezrukavnikov--Finkelberg
\cite{BezFin}; it should be compatible with it, under the action of $\operatorname{Shv}(L^{+}G\backslash LG/L^{+}G)$ on the left, and the moment map $\check{M} \to \check{\mathfrak{g}}^{*}$ on the right. There is a similar, categorical version of the global Conjecture~\ref{conj-global}, for which I defer to the upcoming article.

\section{Beyond endoscopy} \label{sec:BE}

\subsection{Relative functoriality}

\subsubsection{} Let $X, Y$ be two spherical varieties (for possibly different groups $G, G'$), and let $r$ be a morphism of their $L$-groups,
$ r:\LG_X\to \LG_Y.$
According to the relative local Langlands conjecture of \S\ref{conj-local}, it should give rise to a map 
\[ \{X\mbox{-distinguished $L$-packets}\} \longrightarrow \{Y\mbox{-distinguished $L$-packets}\},\]
at least for $L$-packets distinguished in the $L^2$-sense. 

A basic tenet of Langlands' ``beyond endoscopy'' proposal \cite{Langlands-BE}, generalized to the relative setting, states that the resulting map of stable relative characters 
$ J_{\phi_1}^{X} \mapsto J_{\phi_2}^{Y}$
should be realized as the adjoint of a ``transfer operator'' between spaces of stable test measures,
\begin{equation}\label{transfer} \mathcal T: \mathcal S(\mathfrak Y)^\st \to \mathcal S(\mathfrak X)^\st,
\end{equation}
where $\mathfrak Y$ denotes the stack $(Y\times Y)/G'$, and $\mathfrak X = (X\times X)/G$. In most cases, one can take ``stable'' to mean the image of the canonical pushforward map 
\[ \mathcal S(\mathfrak X) \to \operatorname{Measures}((X\times X)\sslash G).\]
When the map $\mathcal S(\mathfrak X)\to \mathcal S(\mathfrak X)^\st$ is an isomorphism (e.g., for the Kuznetsov formula), we will be dropping the exponent ``st.''

In the group case, this operator has been studied by Langlands \cite{Langlands-ST} and Johnstone \cite{Johnstone} when $X$ is a torus and $Y$ is $\GL_n$. 
Understanding these transfer operators could be considered as the basic problem of functoriality, at least in the local setting.

\subsubsection{}  In the global setting, one would have to find a way to employ these transfer operators in a comparison of relative trace formulas. Langlands' proposal, generalized to our setting, is to extract the part of the automorphic spectrum of $\mathfrak Y$ that is is in the image of the functorial lift from $\mathfrak X$ from the (stable) relative trace formula for $\mathfrak Y$  by means of poles of $L$-functions. 

The question of whether it is possible to identify the spectrum of $\mathfrak X$ by orders of poles of $L$-functions has been studied and is known to have a negative answer, in general \cite{AYY}. Other difficulties with this proposal include the isolation of the tempered part of the spectrum; a lot of hard work has gone into this problem, already for the case of $\GL_2$ \cite{FLN, Altug1, Altug2, Altug3}. 

\subsection{An example: symmetric square lift}

\subsubsection{}
Rather than speculating on how to overcome these difficulties, it may be more instructive to look at a variant of the idea, which was applied successfully in the thesis of Venkatesh \cite{Venkatesh}, and to understand what the structure of local transfer operators can tell us about the global problem. Here, $X=T$ is a $1$-dimensional torus over a global field $k$ (the kernel of the norm map for a quadratic etale algebra $E/k$ whose quadratic idele class character we will denote by $\eta$), and $Y$ is the Whittaker model of the group $G=\Gm \times \SL_2$, so that $\check G_Y = \check G$. There is a morphism of $L$-groups
$r: {^LT} \to \LG,$
whose image stabilizes a vector under the product of the standard representation of $\Gm$ with the adjoint representation of $\PGL_2$ (that is, the symmetric-square representation of $\GL_2$, as it factors through $\GL_2 \to \check G = \Gm \times \PGL_2 \to \GL_3$). Let $\mathbb Z/2$ act on $T$ by inversion. The local transfer operator for this morphism was computed in \cite{SaTransfer2}:
\begin{theorem*}
 Let $G, T$ as above be defined over a local field $F$. There is a transfer operator 
 \[\mathcal T:  \mathcal S(N,\psi\backslash G/N,\psi) \to \mathcal S(T)^{\mathbb Z/2},\]
 such that the pullback of every unitary character of $T$ is the Kuznetsov relative character of its functorial lift. In natural coordinates $(r,t)$ for $N\backslash G\sslash N \simeq \Gm \times \mathbb A^1$, it is given by 
 \begin{equation}\label{transfer-torus} (da)^{-1} \mathcal T f(a) =  (dt)^{-1} \lambda(\eta,\psi) \int_r \int_x f\big ( r, \frac{t}{x} \big)  \eta(xrt) \psi(x)\,  dx,
 \end{equation} 
 where $a\in T$ and $t=t(a)$ is its image through an isomorphism $T\sslash(\mathbb Z/2) \simeq \mathbb A^1 \simeq  N\backslash\SL_2\sslash N$, and $\lambda(\eta,\psi)$ is a constant.
\end{theorem*}

What does this theorem tell us about how to extract from the relative trace formula of $Y$ (that is, the Kuznetsov formula of $G$) the part of the spectrum that is due to the torus $T$? Venkatesh \cite{Venkatesh} performs this extraction in two steps, a Poisson summation formula followed by taking the pole of a zeta integral. As explained in \cite[\S10]{SaTransfer2}, the adelic reformulation of the first step is the Poisson summation formula for the Fourier transform corresponding to the inner integral of \eqref{transfer-torus}, while the second step is a global version of the Mellin transform represented by the outer integral.

\subsubsection{}
Thus, we see that understanding the local transfer operators can guide our steps for the global ``beyond endoscopy'' comparisons of trace formulas. Another example of such a comparison is that between the Kuznetsov formula of $G=\GL_2$ and the Selberg trace formula for the same group. The local transfer operator for this comparison was computed in \cite[\S4]{SaTransfer1}, and is given by a simple Fourier tranform (see Theorem \ref{mainthm-G}). Restricted to holomorphic cusp forms, this global comparison via a Poisson summation formula was performed in the thesis of Zeev Rudnick \cite{Rudnick}, about 10 years before Langlands' ``beyond endoscopy'' proposal. A generalization of this comparison to the full Kuznetsov formula, and for $\GL_n$ with $n$ arbitrary, is the object of ongoing joint work with Chen Wan.

\subsubsection{}
It has hopefully become clear that understanding the transfer operators is of paramount importance for the problem of functoriality. In \cite{SaRankone}, I showed that these operators have a very uniform form, for spherical varieties of rank 1. In the remainder of this paper, I would like to propose a reinterpretation of this work, which provides an understanding of those transfer operators as ``change of Schr\"odinger model/geometric quantization'' associated to a symplectic group scheme.

\section{Transfer operators and quantization} \label{sec:transfer} 

The goal of this section is to recast the transfer operators of functoriality, studied in \cite{SaRankone}, in the language of quantization. The idea that quantization should have something to do with functoriality is not new; V.\ Lafforgue suggested it several years ago (private communication), in order to interpret the Rankin--Selberg method, and the functoriality kernels of L.\ Lafforgue \cite{LLafforgue}. Here, however, we apply this idea in a different setting: the setting of ``beyond endoscopy,'' and of the quotient stacks showing up in the relative trace formula -- the hope being that these operators of functoriality will always exist in this setting, even if they do not exist for the spaces ``upstairs.''

Geometric quantization was introduced by Kostant and Souriau \cite{Kostant, Souriau}, following the work of Kirillov on the orbit method \cite{Kirillov}. Since the notion of quantization for measures on stacks that we need has not been developed yet, we will take a phenomenological approach, with ad hoc definitions that provide the desired reformulation of the results of \cite{SaRankone}.

\subsection{Cotangent space of the RTF stack}

\subsubsection{}
The groundbreaking work of Friedrich Knop has shown that, although spherical varieties can be very different from each other, their \emph{cotangent bundles} are quite similar. This will be the basis of our considerations, when we try to relate cotangent bundles of different quotient stacks of the form $(X\times X)/G$. 

For the rest of this paper we will assume, for simplicity, that all groups are split, defined over a field $F$ in characteristic zero. The results of Knop, then, recalled in \S\ref{Knopwork}, hold verbatim over $F$. We assume that $X$ is smooth and quasiaffine, and set again $\mu: M=T^*X \to \mathfrak g^*$ for the cotangent bundle and its moment map.

The ring of regular functions $F[M]$ has a Poisson structure. Knop has studied the subalgebra $F[M]^G$ of $G$-invariants; when $X$ is spherical, this subalgebra is Poisson-commutative, and can be naturally identified with the algebra of regular functions on the affine space $\mathfrak c_X^*$, defined in \S\ref{Knopwork}.
Hence, regular functions on $\mathfrak c_X^*$ pull back to a Poisson-commuting algebra of $G$-invariant functions (``Hamiltonians'') on $M$. One can ask whether the corresponding Hamiltonian vector fields can be integrated to the action of an abelian group scheme $J_X$ (over $\mathfrak c_X^*$) of $G$-automorphisms of $M$, and Knop has answered this in the affirmative \cite{KnAut}.

\subsubsection{}
More precisely, $J_X$ is ``the group scheme of regular centralizers in the split reductive group $G_X$ dual to $\check G_X$.'' When $\check G_X$ is adjoint, so that $G_X$ is simply connected, the group scheme $J_X$ has an explicit description as
\begin{equation}\label{Weilrestriction}
 \Big( \operatorname{Res}_{\mathfrak a_X^*/\mathfrak c_X^*} T^*A_X  \Big)^{W_X} 
\end{equation}
(see \cite[\S2.4]{Ngo-FL}), where $\operatorname{Res}_{\mathfrak a_X^*/\mathfrak c_X^*}$ denotes Weil restriction of scalars from $\mathfrak a_X^*$ to $\mathfrak c_X^*$.  

In general, the group scheme $J_X$ acting on $M$ is an open subgroup scheme of \eqref{Weilrestriction}, which depends not only on the pair $(A_X,W_X)$, but also on the root datum of $X$. Knop defines a slightly different root datum than ours in \cite[\S6]{KnAut}, giving the maximal possible subgroup scheme acting on $M$. For our purposes, we will be content with taking $J_X=$ the open subgroup scheme of \eqref{Weilrestriction} that corresponds to the set of normalized spherical roots, \S\ref{sphericalroots}. This can be described as the regular centralizer group scheme of $G_X$, and 
is the complement of a divisor in \eqref{Weilrestriction}, see \cite[Theorem 7.7]{KnAut}, \cite[\S2.4]{Ngo-FL}.

\subsubsection{} \label{JXrankone}

Let us discuss the rank-$1$ cases. Consider, first, the case $\check G_X = \PGL_2$, so that $A_X=\Gm$, $W_X=\mathbb Z/2$, and the normalized spherical root is twice the generator of the character lattice. (The isomorphism $A_X\simeq \Gm$ is canonical, if we require positive roots to correspond to positive powers.) Then, 
$J_X = J_{\SL_2}$ is given by the restriction of scalars \eqref{Weilrestriction}, which can explicitly be described as follows: Identify $\mathfrak a_X^*=\mathfrak g_m^*$ with the affine line, with coordinate $\sigma = $ the differential of the identity cocharacter, and set $\xi=\sigma^2$, a coordinate on $\mathfrak c_X^*$.  We can write 
\[J_X =  \operatorname{Spec} F[t_0, t_1, \xi] / (t_0^2 -\xi t_1^2 -1), \] 
so that the canonical base change map 
\[ J_X \bullet \mathfrak a_X^* \to T^* A_X = T^*\Gm = \Gm \times \mathfrak g_m^*,\]
where by $\bullet$ we denote fiber product over $\mathfrak c_X^*$,
is given by $(t_0, t_1, \sigma)\mapsto (t_0+\sigma t_1, \sigma)$.

The symplectic form is given by 
\[ \omega = \frac{dt_1 \wedge d\xi}{2 t_0} = \frac{dt_0 \wedge d\xi}{2\xi t_1} = dt_0 \wedge d(t_1^{-1}).  \]
It is immediate to check that this is regular and nondegenerate everywhere on $J_X$.

On the other hand, when the normalized root datum of $X$ is that of $\PGL_2$ (i.e., the normalized spherical root is a generator of the character lattice), the fiber of \eqref{Weilrestriction} over the nilpotent point $0\in \mathfrak c_X^*$ is isomorphic to $\Ga \times \{\pm 1\}$, but the fiber of $J_X$ is just $\Ga$.

\subsubsection{} \label{action} Returning to the general case, the Lie algebra of $J_X$ is canonically isomorphic to the cotangent space of $\mathfrak c_X^*$. Thus, the Hamiltonian vector fields associated to $F[\mathfrak c_X^*]$ give rise to a homomorphism from $\operatorname{Lie}(J_X)$ to $G$-invariant vector fields along the fibers of $M\to \mathfrak g_X^*$ (notation as in \S\ref{invariants}). Knop has shown \cite{KnAut} that these vector fields integrate to an action of $J_X$ on $M=T^*X$ over  $\mathfrak g_X^*$, commuting with the action of $G$. 

Moreover, over a dense open subset $\mathring{\mathfrak c}_X^*\subset \mathfrak c_X^*$, the map $M\to \mathfrak g_X^*$ is a $J_X$-torsor, and the action of $J_X$ arises from the stabilizers of points of $\mathfrak g_X^*$, or even of $\mathfrak g^*$, in $G$, i.e., the stabilizer $G_z$ of a generic point $z\in \mathfrak g^*$ in the image of the moment map acts transitively on the fiber through a map $G_z\to J_X$. 

\subsubsection{}
In order to study the relative trace formula for a stack of the form $\mathfrak X=(X\times X)/G$, I propose to use its cotangent stack
\[ T^*\mathfrak X = (T^*X \times_{\mathfrak g^*} T^*X)/G.\]
(Strictly speaking, the fiber product over $\mathfrak g^*$ should taken with respect to the moment map and its negative, $(\mu, -\mu)$, so that the quotient above corresponds to symplectic reduction with respect to the diagonal action of $G$. We apply multiplication by $-1$ on the fibers of one factor, in order to have fiber product with respect to $(\mu, \mu)$, which is notationally simpler.)

The fiber product, here, should be taken in the derived sense, turning this quotient into a symplectic derived stack. Although the derived features are likely to be important in the future, for the purposes of the current paper we will ignore them. Then, the fiber product over $\mathfrak g^*$ coincides with the fiber product over $\mathfrak g_X^*$ over a dense open $\mathring{\mathfrak c}_X^*\subset \mathfrak c_X^*$, which we will take small enough so that it also has the properties of \S\ref{action}.

Knop's theory, now, gives a very satisfactory description of a dense substack. Namely, consider the diagonal embedding 
$ T^*X \hookrightarrow T^*X \times_{\mathfrak g^*} T^*X,$
which by the action of $J_X$ on the first variable gives rise to 
\begin{equation}\label{JXfiberprod} J_X \bullet T^*X \to T^*X \times_{\mathfrak g^*} T^*X,
\end{equation}
where again by $\bullet$ denotes $\times_{\mathfrak c_X^*}$. Given that $T^*X\to \mathfrak g^*$ is generically a $J_X$-torsor over its image, this map is birational into an irreducible component of the right hand side.

Taking quotients by the $G$-action in \eqref{JXfiberprod}, and using the invariant moment map $T^*X\to \mathfrak c_X^*$, we obtain a correspondence
\begin{equation}\label{fromJX}
 J_X \leftarrow (J_X\bullet T^*X)/G  \to T^*\mathfrak X = (T^*X \times_{\mathfrak g^*} T^*X)/G.
\end{equation}

Over the dense open subset $\mathring{\mathfrak c}_X^*$, the right arrow is an isomorphism, and the left arrow is an isomorphism if we ignore stabilizers. 
We would like to think of \eqref{fromJX} as saying that the symplectic stack $T^*\mathfrak X$ is ``birational'' to the symplectic group scheme $J_X$. Of course, this is a very naive notion of birationality since, even over a good open dense subset, we are ignoring stacky and derived structures in $T^*\mathfrak X$. Nonetheless, even this weak correspondence is quite remarkable, since $J_X$ depends only on the dual group of $X$. \emph{It is not too far-fetched to imagine that this correspondence plays a key role in functoriality.}

\subsection{Rank 1 spherical varieties}

\subsubsection{} We will now specialize to spherical varieties $X = H\backslash G$ with $G$ and $H$ split reductive groups, whose dual group is $\check G_X= \PGL_2$ or $\SL_2$. The group scheme $J_X$ was described in \S\ref{JXrankone}. 

In this setting, \cite{SaRankone} gave explicit formulas for transfer operators \eqref{transfer} between $\mathfrak Y=$ the Kuznetsov stack for the group with dual group $\check G_X$ (see \S\ref{Kuznetsovstack} below), and $\mathfrak X=(X\times X)/G^\diag$. These operators transfer spaces of test measures to each other, but in a number of cases, studied in \cite{SaTransfer1, SaTransfer2, Gan-Wan}, properties such as the transfer of characters or the appropriate fundamental lemma are also known; thus, there is enough evidence to believe that these are the ``correct'' operators of functoriality for these comparisons.

The simplest form of transfer operators appears when the normalized root datum of $X$ is simply connected, i.e., $\check G_X=\PGL_2$. For some varieties, this is possible to achieve by passing to a finite cover; the only cases where this can be done lead to the following spaces:
\begin{equation}\label{typeG}
 X = \SO_{2n-1}\backslash \SO_{2n}, \,\, \Spin_7\backslash \Spin_8, \,\, \operatorname{G}_2\backslash \Spin_7.
\end{equation}
We will call them ``cases of type $G$,'' because the base case is the group variety $\SO_3\backslash \SO_4\simeq \SL_2$. (Here, the $\SO_7\backslash \SO_8\simeq \Spin_7\backslash \Spin_8 \simeq \operatorname{G}_2\backslash \Spin_7$ as varieties, but with the action twisted by the triality automorphism of $\Spin_8$, for the second, and restricted to the subgroup $\Spin_7$, for the third.)

The remaining cases of spherical varieties satisfying our assumptions are 
\begin{equation}\label{typeT}
 X = \GL_n\backslash \PGL_{n+1},\,\, \SO_{2n}\backslash \SO_{2n+1}, \,\, \Sp_{2n-2}\times\Sp_2\backslash \Sp_{2n} \mbox{ (with $n\ge 2$)}, \,\, \Spin_9\backslash F_4, \,\, \SL_3\backslash G_2.
\end{equation}
We will call them ``cases of type $T$,'' because in the base case $\Gm\backslash \PGL_2$ the stabilizer is a torus. (Again, $\SL_3\backslash G_2\simeq \SO_6\backslash \SO_7$, but with the action restricted to $G_2\subset \SO_7$.)

In those cases, the normalized root datum of $X$ is that of $\PGL_2$, and, as we will see, it will be necessary to ``lift'' our description of transfer operators to the root datum of $\GL_2$.

\subsubsection{} \label{Kuzrankone} Let $X$ be as above. Let $G'$ be the split reductive group with the same dual group as $X$, that is, $G'=\SL_2$ for the varieties of \eqref{typeG} and $G'=\PGL_2$ for the varieties of \eqref{typeT}. Let $N\subset G'$ be the upper triangular unipotent subgroup, identified with the additive group $\Ga$, fix a nontrivial character $\psi$ of $F$, and let $Y$ be the Whittaker model of $G'$ with respect to $(N,\psi)$, \S\ref{excellent}. Let $Y^-$ be the Whittaker model with respect to the inverse character, $\psi^{-1}$. 

\label{Kuznetsovstack}
We will symbolically write $\mathfrak Y = (Y\times Y^-)/G'$ for the ``Kuznetsov stack,'' but we will really treat it not as an abstract quotient stack, but as one equipped with the line bundle defined by the Whittaker character. More precisely, this symbol will only find a rigorous meaning in its Schwartz space $\mathcal S(\mathfrak Y)$, which we define to be the $G'^{\diag}$-coinvariant space of the space $\mathcal S(Y\times Y^-)$ of Whittaker Schwartz measures. 

Having identified $N$ with $\Ga$, we let $M'= (f+ \mathfrak n^\perp)\times^N G'$ be the \emph{Whittaker cotangent bundle}, where $f$ is a nilpotent element of $(\mathfrak g')^*$ that is equal to the identity functional on $\mathfrak n = \Ga$. The corresponding bundle for $Y^-$ is $M'^-= (-f+ \mathfrak n^\perp)\times^N G'$. Both come equipped with natural moment maps to $\mathfrak g'^*$, which we will indiscriminately denote by $\mu$. We now define the \emph{Kuznetsov cotangent stack} as 
\[ T^*\mathfrak Y = (M' \times_{\mu, \mathfrak g'^*, (-\mu)} M'^-)/G' \simeq (M' \times_{\mu, \mathfrak g'^*, \mu} M')/G'.\]

Note that the invariant-theoretic quotient $\mathfrak g'\sslash G'$ is canonically identified with $\mathfrak c_X^*$, and the group scheme of regular centralizers in $G'$ is canonically identified with $J_X$. 
It is well-known that $M'$ is a $J_X$-torsor over the regular subset of $\mathfrak g'^*$; in fact, a Kostant section provides a section for this torsor. Therefore, the same considerations that led us to \eqref{fromJX} hold, but here we have an exact isomorphism
\begin{equation}\label{fromJX-Kuz}
J_X \simeq  T^*\mathfrak Y. 
\end{equation}

Our hope, now, is to demonstrate the following idea:
\begin{quote}
 The cotangent stacks $T^*\mathfrak X$ and $T^*\mathfrak Y$ being roughly isomorphic to $J_X$ (by \eqref{fromJX}, \eqref{fromJX-Kuz}), there is a transfer operator of functoriality
 \[\mathcal T: \mathcal S(\mathfrak Y) \to \mathcal S(\mathfrak X)^\st,\]
 corresponding to a ``change of geometric quantization'' for $J_X$.
\end{quote}

\subsubsection{}
Quantization, of course, is as much of a science as an art, and the reader should not expect a rigorous formulation of this hope in this article. In particular, the type of geometric quantization that we need (suitable for encoding measures on stacks) has not, to my knowledge, been developed yet. Therefore, the real content of the results that follow is already contained in \cite{SaRankone}; but we will dress them up in an ad hoc language of quantization, in order to exhibit some deeper structure that seems to be lying behind them.

We will also assume, from now on, that our base field is $F=\mathbb R$, in order to use the language of line bundles with connection, and will write the chosen additive character as $\psi(x) = e^{ i \hbar x}$, where $\hbar$ is a nonzero real constant. The final results, contained in Theorems \ref{mainthm-G}, \ref{mainthm-T}, are valid and were proven in \cite{SaRankone} over an arbitrary local field in characteristic $0$, just by an obvious translation of the formulas. We will only care to describe transfer operators up to an absolute scalar; therefore, we will feel free to choose measures that only modify the result by a scalar, without commenting on those choices.

\subsection{Geometric quantization for type $G$}

\subsubsection{}

The process of geometric quantization on a (real) symplectic manifold $(M,\omega)$ consists in fixing a (complex) Hermitian vector bundle $L$, equipped with a connection $\nabla$ whose curvature is $i \hbar \omega$,  as well as a Lagrangian foliation $\mathscr F$, such that the space of leaves $M/\mathscr F$ is a Hausdorff manifold. 
Then, one attaches to these data the vector space $\mathcal D_{\mathscr F}(M,L)$ of smooth half-densities on $M/\mathscr F$ valued in the space of sections of $L$ over $M$ that are constant along the foliation $\mathscr F$ (with respect to the connection). This space has a canonical inner product (namely, the $L^2$-inner product over $M/\mathscr F$), giving rise to a Hilbert space, by completion. 

To reformulate the results of \cite{SaRankone} in this language, we will now recast the space $\mathcal S(\mathfrak X)^\st$ of stable test measures for the relative trace formula as a space of half-densities on the quotient of the group scheme $J_X$ by a Lagrangian foliation, valued in a line bundle $L_X$.

\subsubsection{} 
Fix a rank-1 space $X$ ``of type $G$,'' i.e., in the list \eqref{typeG}.

Consider the composition of maps
\begin{equation} \label{composition} J_X \bullet T^*X \to T^*X \times_{\mathfrak g^*} T^*X \to X\times X 
\end{equation}
induced from \eqref{JXfiberprod}. There is a natural scaling $\Gm$-action on the left hand side, under which this composition is invariant, and if we consider the ``projectivization'' of the space on the left (= remove the zero section in $T^*X$ and divide by $\Gm$), it was shown in \cite[\S3]{SaRankone} that the resulting map
\begin{equation}\label{resolution} 
\mathbb P(J_X \bullet T^*X) \to X\times X
\end{equation}
is generically an isomorphism. More precisely, in the type-$G$ cases it is an isomorphism over $\mathbb P(J_X^\circ \bullet T^*X)$, where $J_X^\circ\subset J_X$ is the complement of the divisor given by the homogeneous equation $t_1=0$. (This is a combination of Propositions 3.3.2, 3.5.1 in \cite{SaRankone}, and the fact that those spaces have an involutive $G$-automorphism.)

The invariant-theoretic quotient $X\times X \to (X\times X)\sslash G$ is an affine line, and its composition with \eqref{composition} is the map $J_X \bullet T^*X \to J_X \to \mathbb A^1$ that remembers only $t_0$ from the triple $(t_0, t_1, \xi)$ \cite[Proposition 3.4.2]{SaRankone}. We notice that the level sets of $t_0$ on $J_X^\circ$ form a Lagrangian foliation; we will call this foliation ``vertical,'' and denote it by $\mathscr F_\ver$.

\subsubsection{}
Let $d=\dim X$. The short version of the story that follows is that we replace the element $f\in \mathcal S(\mathfrak X)^\st$ (a measure in the variable $t_0$) by 
\begin{equation}\label{meas2dens-G}
 f(t_0)(dt_0)^{-\frac{1}{2}} |t_1|^{-\frac{d}{2}+1}, 
\end{equation}
obtaining a half-density on $J_X^\circ/\mathscr F_\ver$ valued in the line bundle $L_X$ whose sections are functions on $J_X^\circ/\mathscr F_\ver$ multiplied by the factor $|t_1|^{-\frac{d}{2}+1}$. More precisely, $L_X$ will be identified with the trivial line bundle on $J_X^\circ$, but endowed with a connection 
\begin{equation}\label{theconnection}
 \nabla^X = \nabla^0 + d \log|t_1|^{\frac{d}{2}-1} - i\hbar t_1^{-1} dt_0 = \nabla^0 + \Big(\frac{d}{2}-1\Big) t_1^{-1} dt_1 - i\hbar t_1^{-1} dt_0,
\end{equation}
with curvature $i\hbar \omega$, where $\nabla^0$ is the standard flat connection, so that its parallel sections along the vertical foliation are as described. 

Presented this way, this connection is completely unmotivated. In \S\ref{LXnatural} below, we will discuss a more natural description of the pair $(L_X, \nabla^X)$. Continuing, for now, in this ad hoc fashion, \eqref{meas2dens-G} defines a map 
\begin{equation}\label{toLX} 
\mathcal S(\mathfrak X)^\st \to D_\hor (J_X^\circ, L_X),
\end{equation}
where $D_\hor$ (with regular font) denotes \emph{continuous} (not necessarily Schwartz, or even smooth) ``horizontal'' half-densities valued in $L_X$ (i.e., 
half-densities on $J_X^\circ/\mathscr F_\ver$ valued in the descent of $L_X$ by parallel transport). The image of this map will be denoted by $\mathcal D(\mathfrak X)$.

\subsubsection{} \label{Ytodensities}
We consider another Lagrangian foliation $\mathscr F_\hor$ on $J_X^\circ$, which we will call ``horizontal:'' its leaves are the level sets of $t_1$. 
For the line bundle with connection $(L_X, \nabla^X)$, as above, flat sections along horizontal leaves are simply functions of $t_1\ne 0$, multiplied by the factor $\psi(\frac{t_0}{t_1})$; note that this description is independent of the dimension $d$ used to define $\nabla^X$.

We now propose to think of the space of test measures $\mathcal S(\mathfrak Y)$ for the Kuznetsov formula as a subspace $\mathcal D(\mathfrak Y)\subset D_\ver(J_X^\circ, L_X)$, where $D_\ver$ denotes continuous ``vertical'' half-densities (i.e., half-densities on $J_X^\circ/\mathscr F_\hor$) valued in $L_X$.
First of all, consider the map 
\[J_X \bullet M' \xrightarrow\sim M'\times_{\mathfrak g'^*} M',\]
with notation as in \S\ref{Kuzrankone}, where the action of $J_X$ is again on the first copy of $M'$.

\begin{lemma*} \label{Ycoordinates}
The composition of the map above with $M'\times_{\mathfrak g'^*} M' \to Y\times Y \to (Y\times Y)\sslash G$ is the map that only remembers the coordinate $t_1$ of $J_X$.
\end{lemma*}
 This is the reason why the foliation $\mathscr F_\hor$ is relevant to the Kuznetsov formula. \withproofs{ The proof, and similar calculations later, will use the identification $\mathfrak{sl}_2^*\simeq \mathfrak{sl}_2$ through the Killing form $\langle A, B\rangle \mapsto \operatorname{tr}(AB)$. }

\begin{proof}
  Let $\xi=\sigma^2 \in \mathfrak c_X^*$ correspond to a split regular semisimple conjugacy class, represented by the element $\begin{pmatrix} \frac{\sigma}{2} \\ 1 & -\frac{\sigma}{2}\end{pmatrix} \in \mathfrak g'^*$, which belongs to the moment image of the fiber of $T^*Y$ over the point $N1\in Y$; we will denote the corresponding point on $T^*_{N1}Y$ by $\tilde \xi$. The fiber of $J_X$ over $\xi$ acts on $\tilde \xi \in T^*Y$ via the centralizer of this element, which contains elements of the form $\begin{pmatrix} t_0 + \sigma t_1 \\ 2 t_1 & t_0-\sigma t_1\end{pmatrix}$ (with $(t_0,t_1, \xi=\sigma^2)$ as before). We can parametrize $(Y\times Y)\sslash G' = N\backslash G' \sslash N$ by the variable $\begin{pmatrix} a & b \\ c & d \end{pmatrix} \mapsto c$; thus, acting by $J_{X,\xi}$ on the first copy of $\tilde\xi$ in $(\tilde\xi, \tilde\xi)\in M'\times M'$, we arrive at a point with $c = 2 t_1$. 
\end{proof}

There is a natural pullback from Whittaker functions on $Y\times Y^-$ to scalar-valued functions on $M'\times_{\mathfrak g'^*} M'$, as follows:
Thinking of elements of $\mathcal F(Y)$ (that is, Whittaker functions) as sections of a line bundle $L_\psi$ over $Y= N\backslash G'$ (and similarly for $Y^-$, just replacing $\psi$ by $\psi^{-1}$), we note that the line bundle $L_\psi \boxtimes L_{\psi^{-1}}$ is canonically trivial over the diagonal $Y^\diag\subset Y\times Y^-$. There is now a unique trivialization of its pullback to $J_X \bullet M'$ 
which coincides with the canonical one over the diagonal, and is equivariant with respect to the action of $J_X\times G'$. More explicitly, if we use a Kostant section to identify $T^*Y \simeq \mathfrak c_X^*\times G'$, and the negative of that section for $Y^-$, we pull back Whittaker functions to scalar-valued functions on $T^*Y$ via the projection to $G'$, and then restrict to $T^*Y \times_{\mu, \mathfrak g'^*, (-\mu)} T^*Y^- \simeq M' \times_{\mathfrak g'^*} M'$.

The short version of the story, now, is that we fix a $G'$-invariant measure on $N\backslash G'$, use it to identify Schwartz (Whittaker) measures on $Y\times Y^-$ with Schwartz (Whittaker) functions, pull them back to scalar-valued functions on $J_X^\circ \bullet M' \simeq J_X^\circ \times G'$, and integrate them against a chosen Haar measure on $G'$. This gives functions on $J_X^\circ$ that, as can be easily confirmed, correspond to sections of $L_X$, flat along the leaves of $\mathscr F_\hor$; further multiplying them by the factor $|t_1|^\frac{1}{2} dt_1^\frac{1}{2}$ gives rise to an element of 
$D_\ver(J_X^\circ, L_X)$. This descends to an injective map 
\begin{equation}\label{toLY} 
\mathcal S(\mathfrak Y)\to D_\ver(J_X^\circ, L_X),
\end{equation}
whose image will be denoted by $\mathcal D(\mathfrak Y)$. Again, the factor $|t_1|^\frac{1}{2} dt_1^\frac{1}{2}$ seems unmotivated, and we will attempt to explain it in \S\ref{LYnatural} below, after formulating the main theorem.

\subsubsection{} \label{mainthm-G}
To recap, we have defined a line bundle $L_X$ on $J_X^\circ$, endowed with a connection $\nabla^X$ with curvature $i\hbar \omega$,  ``vertical'' and  ``horizontal'' foliations $\mathscr F_\ver$, $\mathscr F_\hor$ on $J_X^\circ$,  and have identified the spaces $\mathcal S(\mathfrak X)^\st$, $\mathcal S(\mathfrak Y)$ of test measures for the corresponding quotients with spaces 
$\mathcal D(\mathfrak X)$, $\mathfrak D(\mathfrak Y)$ of ``horizontal'' and ``vertical'' half-densities for $(L_X,J_X^\circ)$.
The main result \cite[Theorem 1.3.1]{SaRankone} for the transfer operator in this case can now be formulated as follows:

\begin{theorem*}
 There is an injective operator
 \[ \mathcal T: \mathcal D(\mathfrak Y) \to \mathcal D(\mathfrak X)\]
 given by integration along the leaves of the vertical foliation: 
 \[D_\ver (J_X^\circ, L_X) \dashedrightarrow D_\hor (J_X^\circ, L_X),
\]
\begin{equation}\label{transferG} \mathcal T\varphi (j) = \int_{\mathscr F_{\ver,j}} T_{j,j'}(\varphi(j')) |\omega(j')|^{\frac{1}{2}} ,
\end{equation}
where $\mathscr F_{\ver,j}$ denotes the leaf of $\mathscr F_\hor $ through the point $j$, and $T_{j,j'}$ denotes parallel transport from the fiber of $L_X$ over $j'$ to the fiber over $j'$ along this leaf.

Its inverse $\mathcal T^{-1}$, valued in an enlargement $\mathcal D_{L_X}^-(\mathfrak Y)\supset \mathcal D(\mathfrak Y)$ described in \cite[\S1.3]{SaRankone}, is given by integration along the leaves of the horizontal foliation:
\[D_\hor (J_X^\circ, L_X) \dashedrightarrow D_\ver (J_X^\circ, L_X),
\]
\end{theorem*}

Note that a horizontal half-density on $J_X^\circ$, multiplied by the half-density $|\omega|^\frac{1}{2}$, gives rise to a vertical half-density valued in the bundle of densities on the leaves of $\mathscr F_\hor$; thus, it makes sense to integrate it along these leaves, obtaining a half-density on $J_X^\circ/\mathscr F_\hor$. This, of course, is completely analogous to the canonical intertwiners for the Schr\"odinger models quantizing a symplectic vector space \cite{WWLi}.

\begin{proof}
 Let us start with an element $h\in \mathcal S(\mathfrak Y)$. In \cite{SaRankone}, those were considered as measures in the variable $2t_1$, by trivializing them along the representatives $\widetilde{t_1}:= \begin{pmatrix} & -(2t_1)^{-1} \\ 2t_1 \end{pmatrix}$ for $N\backslash G'/N$-cosets. Let us write this measure as $\Phi(\widetilde{t_1}) \cdot  |t_1| dt_1$, and extend the function $\Phi$ to the unique function on $N\widetilde{t_1}N$ which satisfies $\Phi(n_1 \widetilde{t_1} n_2) = \psi(n_1 n_2) \Phi(\widetilde{t_1})$.  
 The map \eqref{toLY} translates this to the vertical half-density $\varphi(j) = \tilde\Phi(j) |t_1|^\frac{1}{2} dt_1^{\frac{1}{2}}$ on $J_X^\circ$, where $\tilde\Phi$ is the pullback of $\Phi$ to $J_X^\circ$ through the map 
 \begin{equation}\label{fromJXKostant}
J_X \to T^*Y \simeq \mathfrak c_X^*\times G'\to G'  
 \end{equation}
determined by a Kostant section.

To compute its integral along the leaf $\mathscr F_{\ver,t_0}$ of the vertical foliation with coordinate $t_0$, recall that parallel translation along the vertical leaves with respect to the connection $\nabla^X$ (starting from a point with $t_1=1$) translates to multiplication by $|t_1|^{-\frac{d}{2}+1}$ in this trivialization. Therefore, using the fact that $|\omega|^\frac{1}{2} = |t_1|^{-1} |dt_1\wedge dt_0|^\frac{1}{2}$, the scalar-valued expression for this integral is
 \begin{align*} \mathcal T\varphi(t_0,t_1,\xi) &= \underset{\nabla^X}{\underbrace{|t_1|^{-\frac{d}{2}+1}}} \int_{\mathscr F_{\ver,t_0}}  \underset{\nabla^X}{\underbrace{|t_1'|^{\frac{d}{2}-1}}}\cdot \underset{\mbox{\small triv.\ for }Y}{\underbrace{\tilde\Phi(t_0,t_1',\xi') |t_1'|^\frac{1}{2} (dt_1')^{\frac{1}{2}}}} \cdot \underset{|\omega|^\frac{1}{2}}{\underbrace{|t_1'|^{-1} (dt_1')^{\frac{1}{2}}\cdot (dt_0)^\frac{1}{2}}} \\ 
 &=   (dt_0)^\frac{1}{2}\cdot |t_1|^{-\frac{d}{2}+1} \int_{\mathscr F_{\ver,t_0}} \tilde\Phi(t_0,t_1',\xi') |t_1'|^\frac{d-3}{2} \, dt_1'.
\end{align*}
Using the identification $\mathfrak{sl}_2^*\simeq \mathfrak{sl}_2$ via the Killing form, and the Kostant section $\xi\mapsto \begin{pmatrix} & \frac{\xi}{4} \\ 1 \end{pmatrix}$, one computes that the map \eqref{fromJXKostant} is 
\[ (t_0,t_1,\xi)\mapsto \begin{pmatrix} t_0 & \frac{\xi}{4}\cdot 2t_1 \\ 2t_1 & t_0\end{pmatrix} = \begin{pmatrix} 1 & \frac{t_0}{2t_1} \\ & 1 \end{pmatrix} \begin{pmatrix} & -(2t_1)^{-1} \\ 2t_1  \end{pmatrix}\begin{pmatrix} 1 & \frac{t_0}{2t_1} \\ & 1 \end{pmatrix} ,\]
hence the expression above reads 
\[(dt_0)^\frac{1}{2}\cdot |t_1|^{-\frac{d}{2}+1} \int_{\mathscr F_{\ver,t_0}} \Phi\begin{pmatrix} & -(2t_1)^{-1} \\ 2t_1  \end{pmatrix} \psi\Big(\frac{t_0}{t_1'}\Big) |t_1'|^\frac{d-3}{2} \, dt_1'.\]
By \eqref{meas2dens-G}, this translates to the measure 
\[dt_0 \cdot \int_{\mathscr F_{\ver,t_0}} \Phi\begin{pmatrix} & -(2t_1)^{-1} \\ 2t_1  \end{pmatrix} \psi\Big(\frac{t_0}{t_1'}\Big) |t_1'|^\frac{d-3}{2} \, dt_1'\]
in the variable $dt_0$. 
This is the formula \cite[(10)]{SaRankone} for the transfer operator -- note that the coordinates used there are $2t_1$ for $N\backslash G'\sslash N$ and $2t_0$ for $(X\times X)\sslash G$. I leave the verification for the inverse operator to the reader.
\end{proof}

\subsubsection{} \label{LXnatural}

The line bundle $L_X$, with its connection, admits a more natural description as the dual to \emph{a line bundle of half-densities on the fibers of the invariant moment map} $\mu_G:M=T^*X\to \mathfrak c_X^*$. The map \eqref{toLX}, then, admits a more natural description as descending, up to a choice of invariant measure on $X$, from a map from Schwartz half-densities on $X\times X$,
\begin{equation}\label{densityquotient}
  \mathcal D(X\times X)\to D_\hor (J_X^\circ, L_X).
\end{equation}
Let us see how this works.

It will be convenient to choose a section $s$ of the invariant moment map $\mu_G$. Such a section exists in the cases $X=H\backslash G$ of \eqref{typeG} when $G$ and $H$ are split; it suffices to check the case of $\SO_{2n-1}\backslash \SO_{2n}$, and we refrain from attempting to give an abstract argument. See \S\ref{nonsplit} for a further discussion of this issue.

Obviously, the section $s$ has image in the smooth locus of the map $\mu_G$, which implies that the fibers of this map are transversal to the section. If $\mathcal O_\xi$ denotes the fiber over $\xi\in \mathfrak c_X^*$, let $D_X$ be the algebraic line bundle over $\mathfrak c_X^*$ whose fiber over $\xi$ is the determinant of the tangent space of $\mathcal O_\xi$ at $s(\xi)$. Let $L_X = |D_X|^\frac{1}{2}$, a complex line bundle whose fiber over $\xi$ is dual to the space of Haar half-densities on this tangent space. By pullback, we will also consider $L_X$ as a line bundle over $J_X$.

\label{LXtrivialization}
For $\xi\ne 0$, the fibers of the invariant moment map $\mu_G$ are $G$-orbits, therefore the tangent space of $\mathcal O_\xi$ is $\mathfrak g/ \mathfrak g_{s(\xi)}$, and the fiber of $L_X$ over $\xi$ is the complex line $|\det \mathfrak g \otimes \det \mathfrak g^*_{s(\xi)}|^\frac{1}{2}$. Its dual is the line of invariant half-densities on $\mathcal O_\xi$.

There is a natural way to trivialize the bundle $L_X$, up to a scalar. It uses the fact that the stabilizers $G_{s(\xi)}$, for $\xi\ne 0$, are isomorphic over the algebraic closure. Thus, $G$-conjugacy gives canonical isomorphisms between the complex line bundles $|\det  \mathfrak g_{s(\xi)}|$, and this allows us to uniformly fix an invariant measure $d\dot g$ on all orbits $\mathcal O_\xi$, for $\xi\ne 0$, see \cite[\S4]{SaRankone}. Then, by \cite[Theorem 4.0.3]{SaRankone}:

\begin{proposition*}
For a suitable choice of $d \dot g$, as above, and the canonical measure $dz$ on $T^*X$ induced by the symplectic form, we have the integration formula
\[ \int_{T^*X} \Phi(z) dz = \int_{\mathfrak c_X^*} |\xi|^{\frac{d}{2}-1} \int_{\mathcal O_\xi} \Phi(s(\xi) \dot g) \, d\dot g \, d\xi.\]
\end{proposition*}

Hence, the family of Haar half-densities $\xi\mapsto \Big (|\xi|^{\frac{d}{2}-1}  d\dot g \Big )^\frac{1}{2}$ on the orbits $\mathcal O_\xi$, for $\xi\ne 0$, extends to a nonvanishing half-density on the fiber over $0$. 
We now use this family (depending up to a constant on our choice of $d\dot g$) to trivialize $L_X^*$, hence also $L_X$, i.e., we have an isomorphism 
\begin{equation}\label{trivial-X}
 L_X \simeq \underline{\mathbb C}
\end{equation}
with the trivial line bundle. Moreover, the proposition above shows that a nonzero element of the fiber of $L_X^*$ over $0$ corresponds to a unique half-density on $\mathcal O_0$, obtained as the limit of $G$-invariant half-densities over the fibers $\mathcal O_\xi$ with $\xi \ne 0$. Hence, each element in the total space of $L_X^*$ gives rise to a half-density on the corresponding fiber of $\mu_G$.

\subsubsection{} 
We can now define the map \eqref{densityquotient}. 
Let $\varphi\in \mathcal D(X\times X)$. The product $\varphi\cdot (dt_0)^{-\frac{1}{2}}$ restricts to a half-density on each fiber of the smooth locus of the invariant-theoretic quotient $X\times X \to  (X\times X)\sslash G$. The idea is to integrate this half-density, but for that purpose we need to turn it into a measure. We will do so after pulling it back to $J_X^\circ$ via the maps
\[ J_X^\circ \bullet (T^*X\smallsetminus X) \to \mathbb P(J_X^\circ \bullet T^*X) \to X\times X,\]
where the second arrow is \eqref{resolution}, an isomorphism onto its image. 

For every $j\in J_X^\circ$ with image $\xi(j) \in \mathfrak c_X^*$, the map \eqref{composition} restricts to a map $\{j\} \bullet \mathcal O_{\xi(j)} \to X\times X$ that is, by \eqref{resolution}, an isomorphism onto its image (up to removing the zero section of $T^*X$, if $\xi(j)=0$). Thus, the pullback of $\varphi\cdot (dt_0)^{-\frac{1}{2}}$ induces a half-density on $\mathcal O_{\xi(j)}$. Multiplying by the half-density corresponding to an element of the fiber of $L_X^*$ over $\xi(j)$ gives rise to a measure, which we can integrate. This way, we get a canonical map
\[ \mathcal D(X\times X) (dt_0)^{-\frac{1}{2}} \times \underline{L_X^*} \to \mathbb C,\]
where $\underline{L_X^*}$ denotes the total space of the line bundle $L_X^*$ over $J_X^\circ$. This corresponds to a map 
\begin{equation}\label{tosections}
\mathcal D(X\times X) (dt_0)^{-\frac{1}{2}} \to  \Gamma(J_X^\circ, L_X),
\end{equation}
where the right hand side denotes (continuous) sections of $L_X$ over $J_X^\circ$. 

\subsubsection{} 
The image of $\varphi\cdot (dt_0)^{-\frac{1}{2}}$ under this map has an invariance property: namely, its values at different points $j$ with the same value of $t_0$ ``coincide.''  To make sense of this, we need to endow $L_X$ with a connection, whose parallel sections along vertical Lagrangians descend to duals of half-densities on the corresponding $G$-orbits on $X\times X$. To explicate this, consider the integration formula of \cite[Theorem 4.0.2]{SaRankone}:

\begin{proposition*}
 For a function $\Phi$ and an invariant measure $dx$ on $X\times X$, we have 
\begin{equation}\label{integration-G} \int_{X\times X} \Phi(x)\, dx = \int_{(X\times X)\sslash G} |t_0^2-1|^{\frac{d}{2}-1} \int_{\mathcal O_{t_0}} \Phi(\dot g)  \, d\dot g \, dt_0.
\end{equation}
\end{proposition*}
Here, we have denoted by $\mathcal O_{t_0}$ the preimage of $t_0 \in (X\times X)\sslash G$ in $X\times X$, using similar notation as for the preimages of points of $\mathfrak c_X^*$ in $T^*X$. The reason is that, as above, for $j\in J_X^\circ$ with $\xi(j)\ne 0$, we can identify the $G$-orbit $\mathcal O_{\xi(j)}$ with the image of $\{j\}\bullet \mathcal O_{\xi(j)}$ in $X\times X$, which is equal to $\mathcal O_{t_0(j)}$. The measure $d\dot g$ in the Proposition, then, is the same measure on $\mathcal O_{\xi(j)}$ as the one used to trivialize the bundle $L_X$ in \S\ref{LXtrivialization}.

The proposition above tells us that, if for every $j$ with $\xi(j)\ne 0$ we multiply $\varphi (dt_0)^{-\frac{1}{2}}$ by the half-density $\Big(|t_0^2(j)-1|^{\frac{d}{2}-1} d\dot g\Big)^\frac{1}{2}$, the integral will depend only on the function $t_0(j)$ of $j$. On the other hand, the trivialization \eqref{trivial-X} of $L_X$ uses the half-density $\Big(|\xi|^{\frac{d}{2}-1} d\dot g\Big)^\frac{1}{2}$. The quotient of the two is $|t_1|^{\frac{d}{2}-1}$. We conclude that the map \eqref{tosections}, composed with the trivialization \eqref{trivial-X}, gives rise to \emph{functions $f$ on $J_X^\circ$ such that $|t_1|^{\frac{d}{2}-1} f$ is constant along fibers of $t_0$}. This explains the definition of $\nabla^X$ in \eqref{theconnection}, and  completes the construction of the map \eqref{densityquotient}.

\subsubsection{} \label{LYnatural}

In a similar way, we define a line bundle $L_Y$ on $J_X$, pulled back from $\mathfrak c_X^*$, as the dual of the line bundle of $G'$-invariant half-densities on the fibers of $M'=T^*Y\to \mathfrak c_X^*$. Here, the fibers are $G'$-torsors, hence fixing a Haar half-density on $G'$ gives rise to a trivialization 
\begin{equation}\label{trivial-Y}
 L_Y \xrightarrow\sim \underline{\mathbb C}.
\end{equation}
Through the trivializations \eqref{trivial-X}, \eqref{trivial-Y}, the line bundles $L_X$ and $L_Y$ are identified. 

Recall from \S\ref{Ytodensities} that the description of ``vertical'' half-densities for $(L_X,\nabla^X)$ is the same in every case (does not depend on the dimension $d$ of $X$). We can now define a map 
\begin{equation}\label{densityquotient-Y}
  \mathcal D(Y\times Y^-) \to D_\ver(J_Y^\circ, L_Y),
\end{equation}
in a completely analogous way to \eqref{densityquotient}, using also the scalar-valued pullback from Whittaker functions on $Y\times Y^-$ to scalar-valued functions on $M'\times_{\mathfrak g'^*} M'$, described in \S\ref{Ytodensities}. To describe it explicitly, 
consider the integration formula for functions on $Y^2$, analogous to \eqref{integration-G},
\begin{equation}\label{integration-Y}
 \int_{Y\times Y} |\Phi(y)| \, dy = \int_{(Y\times Y)\sslash G'} |t_1| \int_{G'} |\Phi(\widetilde{t_1} g)| \, dg  \, dt_1
\end{equation}
(where $\widetilde{t_1}$ denotes any lift of $t_1$ to $Y\times Y$). Symbolically, we can write the $G'\times G'$-invariant measure $dy$ on $Y\times Y$ as $dg \times |t_1| dt_1$. Similarly, we can write the spaces of Schwartz measures and half-densities as 
\begin{align*}
 \mathcal S(Y\times Y) &= \mathcal F(Y\times Y) \cdot \big(dg \times |t_1| dt_1\big),\\
  \mathcal D(Y\times Y) &= \mathcal F(Y\times Y) \cdot \big(dg \times |t_1| dt_1\big)^\frac{1}{2}.
\end{align*}
Choosing the half-density $(dg)^\frac{1}{2}$ to define the trivialization \eqref{trivial-Y}, the image of a Schwartz half-density $\Phi \cdot \big(dg \times |t_1| dt_1\big)^\frac{1}{2}$ under \eqref{densityquotient-Y} is precisely the image of the Schwartz measure $\Phi \cdot \big(dg \times |t_1| dt_1\big)$ described in the definition of \eqref{toLY}, giving a natural meaning to that definition.

\subsection{Geometric quantization for type $T$} \label{sstypeT}

\subsubsection{} 
In the cases \eqref{typeT} of ``type $T$,'' the analog of Theorem \ref{mainthm-G} does not directly hold. In turns out, however, that there is a similar interpretation of transfer operators, if we pass from $J_X$, the group scheme of regular centralizers in $\PGL_2$, to $\widetilde{J_X}=$ the group scheme of regular centralizers in $\GL_2$. It lives over $\widetilde{\mathfrak c_X}^*:= \mathfrak{gl}_2^*\sslash \GL_2$. 

\label{coordinatesGL2} If we write $\GL_2 = (\SL_2 \times \Gm)/\mu_2$, use coordinates $(t_0, t_1, \xi)$, as before, for $J_{\SL_2}$, and coordinates $(z,\tau)$ for $T^*\Gm = \Gm \times \mathfrak g_m^*$, we obtain
\[\widetilde{J_X} =  \operatorname{Spec} F[t_0, t_1, \xi, z^{\pm 1}, \tau]^{\mu_2} / (t_0^2 -\xi t_1^2 -1), \] 
where $-1\in \mu_2$ acts by $(t_0, t_1, \xi, z, \tau)\mapsto (-t_0, -t_1, \xi, -z, \tau)$. The map to $J_X = J_{\PGL_2}$ is then obtained by symplectic reduction modulo $\Gm$, and the symplectic form on $\widetilde{J_X}$ reads 
\[ \omega = dt_0 \wedge d(t_1^{-1}) + d^\times z \wedge d\tau,\]
where the notation is $d^\times z := d \log z = \frac{dz}{z}$.

\subsubsection{} \label{GmGL2}
To motivate the passage to $\widetilde{J_X}$, we should first look at the simple case $X= \Gm\backslash \PGL_2$. This space is a quotient of $\tilde X= \Gm\backslash \GL_2$, where $\Gm$ is embedded as the general linear group of a 1-dimensional subspace, and one can think of $\mathcal S((X\times X)/G)$ as the $\Gm$-coinvariants of the space $\mathcal S((\tilde X\times \tilde X)/\tilde G)$, where $\Gm$ stands for the center of $\tilde G= \GL_2$,  
\[ \mathcal S((X\times X)/G) = \mathcal S((\tilde X\times \tilde X)/\tilde G)_{\Gm}.\]

In this setting, one can easily study a transfer operator 
\[ \mathcal T: \mathcal S(\tilde{\mathfrak Y})\to \mathcal S(\tilde{\mathfrak X}),\]
where $\tilde{\mathfrak X} = (\tilde X\times \tilde X)/\tilde G$, and $\tilde{\mathfrak Y}$ is the ``Kuznetsov quotient stack'' of $\GL_2$, via the ``unfolding'' method.
The ``unfolding'' method \cite[\S9.5]{SV} gives rise to an explicit morphism of Schwartz half-densities
\begin{equation}\label{unfolding} 
\mathcal U: \mathcal D(\tilde Y) \to  \mathcal D(\tilde X)
\end{equation}
(where $\tilde Y$ denotes the Whittaker model of $\GL_2$), which extends to an $L^2$-isometry. 

We can repeat the earlier constructions to identify the spaces of test measures above with spaces of ``vertical,'' resp.\ ``horizontal'' half densities on $\widetilde{J_X}^\circ = $ the complement of $t_1=0$, valued in a line bundle $L_X$ (with connection). Here, the corresponding foliations are determined by the following:

\begin{lemma*}
 The invariant-theoretic quotient $(\tilde X\times \tilde X)\sslash \tilde G= \Gm\backslash \GL_2\sslash \Gm$ is a two-dimensional affine space. One can choose the coordinates $(x,y)$ so that the resulting map 
\[ \widetilde{J_X} \to \widetilde{J_X} \bullet T^*\tilde X \to T^*\tilde X \times_{\tilde{\mathfrak g}^*} T^*\tilde X \to (\tilde X\times \tilde X)\sslash \tilde G\]
(the first arrow again by choosing a section of $T^*\tilde X \to \widetilde{\mathfrak c_X}^*$) is given by 
\begin{equation}\label{2dim} 
 (t_0, t_1, \xi, z, \tau) \mapsto \Big ( x= z^{-1}(t_0 - \tau t_1) , y = z(t_0 + \tau t_1) \Big ).
\end{equation}

The invariant-theoretic quotient for the Kuznetsov formula of $\GL_2$, composed with the analogous map from $\widetilde{J_X}$, 
\[ \widetilde{J_X} \to \widetilde{J_X} \bullet T^*\tilde Y \to N\backslash \tilde G\sslash N\]
is the map that remembers all even-order monomials in the coordinates $t_1$ and $z^{\pm 1}$.

\end{lemma*}

\begin{proof}
 Identify $\mathfrak{gl}_2^*\simeq \mathfrak{gl}_2$ via the trace pairing $\langle A, B \rangle = \operatorname{tr}(AB)$, and use the section $(\xi,\tau)\mapsto \widetilde{(\xi,\tau)} :=\begin{pmatrix} 0 & \frac{\xi-\tau^2}{4} \\ 1 & \tau \end{pmatrix}$ of the map $\tilde{\mathfrak g}^* \to \widetilde{\mathfrak c_X}^*$, to lift $\widetilde{\mathfrak c_X}^*$  to $\mathfrak g_m^\perp = $ the fiber of $T^*\tilde X$ over the point $\Gm 1$ (where $\Gm$ is embedded as the top left entry in $\GL_2$). It is enough to verify the lemma for split semisimple classes $\xi = r^2$, $r\ne 0$, in which case the centralizer of $\widetilde{(\xi,\tau)}$ contains an element of the form
 $\begin{pmatrix} 
z (t_0 - \tau t_1) & * \\
* & z(t_0 + \tau t_1)                                                                                                                                                                                                                                                                                                                                                                                                                                                                                                                                                                                                                                                                                  \end{pmatrix}$, with eigenvalues $(z (t_0 + rt_1), z(t_0 - rt_1))$. Choosing the coordinates $ \begin{pmatrix} A & B \\ C & D \end{pmatrix} \mapsto (x=A \det^{-1}, y = D)$ for $(\tilde X\times \tilde X)\sslash \tilde G= \Gm\backslash \GL_2\sslash \Gm$, this element has $(x, y) = \Big ( z^{-1}(t_0 - \tau t_1) ,  z(t_0 + \tau t_1) \Big)$. 

The calculation for the Kuznetsov quotient follows from Lemma \ref{Ycoordinates}, and the isomorphism $\GL_2 \simeq (\SL_2\times \Gm)/\mu_2$.
\end{proof}

Note that the map $\widetilde{J_X} \to (\tilde X\times \tilde X)\sslash \tilde G$ is smooth, when restricted to $\widetilde{J_X}^\circ= $ the complement of the divisor $t_1=0$. The ``vertical'' foliation $\mathscr F_\ver$ is defined as the set of fibers of this map. The ``horizontal'' foliation $\mathscr F_\hor$ on $\widetilde{J_X}^\circ$ is defined by the level sets of $(t_1 z, z^2)$.

\begin{remark*}
The passage to $(X\times X)\sslash G= \mathbb A^1$ is given by the coordinate
\begin{equation}\label{1dim} 
 (x, y) \mapsto c:= xy = t_0^2 - \tau^2 t_1^2 = (\xi - \tau^2) t_1^2 +1.
\end{equation}
\end{remark*}

\subsubsection{}
Still in the case of $\tilde X= \Gm\backslash \GL_2$, defining line bundles $L_X$, $L_Y$ over $\widetilde{J_X}^\circ $ exactly as before, 
we can repeat the constructions of the maps \eqref{densityquotient}, \eqref{densityquotient-Y} for $\tilde G$, to identify test measures as spaces 
\[\mathcal D(\tilde{\mathfrak X}) \hookrightarrow D_\hor (\widetilde{J_X}^\circ, L_X)\] 
\[\mathcal D(\tilde{\mathfrak Y}) \hookrightarrow D_\ver (\widetilde{J_X}^\circ, L_Y) \] 
of ``horizontal,'' resp.\ ``vertical'' half densities valued in those line bundles. 

\label{mainthm-T}
With the appropriate identification $L_X \simeq L_Y$ over $\widetilde{J_X}^\circ$ (which we will present for the general case in \S\ref{trivial-X2}), we can now descend the unfolding map \eqref{unfolding}, applied to $\mathcal D(\tilde Y)\otimes \mathcal D(\tilde Y^-)$, to coinvariants for the diagonal action of $\tilde G$, obtaining a transfer operator, which can be explicitly described, along the lines of \cite[Theorem 5.4]{SaBE1}:
\begin{theorem*} \label{theorem2d}
 The transfer operator 
\[\tilde{\mathcal T}: \mathcal D(\tilde{\mathfrak Y}) \xrightarrow\sim \mathcal D(\tilde{\mathfrak X})\]
 is the operator of integration along the leaves of the vertical foliation:
\[D_\ver (\widetilde{J_X}^\circ, L_X) \dashedrightarrow  D_\hor (\widetilde{J_X}^\circ, L_X).
\]

The transfer operator
 \[ \mathcal T: \mathcal D(\mathfrak Y) \to \mathcal D(\mathfrak X)\]
is the descent of $\tilde{\mathcal T}$ to $\Gm$-coinvariants. 
\end{theorem*}

 \withproofs{ Since this is included only for motivational purposes, the proof is left to the reader. }

\subsubsection{} \label{trivial-X2}
Let, now, $X$ be a general space from the list \eqref{typeT}. The idea is to generalize the statement of Theorem \ref{theorem2d} for the transfer operator $\mathcal T$, \emph{even though the space $X$ does not, in general, admit a cover such as $\tilde X$}. (Such a cover exists, more generally, for $X=\GL_n\backslash \PGL_{n+1}$, and can be used to motivate some of the definitions that follow.) 

An important feature of the general case is that we will extend the map $J_X\to (X\times X)\sslash G \simeq \mathbb A^1$ to a map $\widetilde{J_X}\to \mathbb A^1$, given by the same coordinate $c$, \eqref{1dim}, as in the case of $\Gm\backslash\GL_2$, 
and will define the vertical and horizontal foliations on $\widetilde{J_X}^\circ$ as in \S\ref{GmGL2}, e.g., the vertical foliation consists of level sets of the pair of functions $(x,y)$ of \eqref{2dim}. We define the complex line bundle $L_X$ over $\mathfrak c_X^*$ as in \S\ref{LXnatural}, and \emph{extend it to a line bundle $L_X$ on $\widetilde{\mathfrak c_X}^*$},  by pullback along the map 
\begin{equation}\label{ctildetoc}\widetilde{\mathfrak c_X}^*\ni (\xi, \tau) \mapsto \xi-\tau^2\in \mathfrak c_X^*.
\end{equation}
 By pullback from $\widetilde{\mathfrak c_X}^*$, this also becomes a line bundle over $\widetilde{J_X}$. We endow it with the same trivialization \eqref{trivial-X} as before. 

Roughly speaking, now, if we fix a Haar measure on $F^\times$, the map \eqref{ctildetoc} allows us to pull back elements of $D_\hor(J_X^\circ, L_X)$ to $\Gm$-invariant elements of $D_\hor(\widetilde{J_X}^\circ, L_X)$, thus obtaining maps
\begin{equation}\label{lifttotilde}
 \mathcal D(X\times X) \to D_\hor(J_X^\circ, L_X) \to D_\hor(\widetilde{J_X}^\circ, L_X)^{\Gm}.
\end{equation}
(Fixing a Haar measure on $F^\times$ allows the switch from the coinvariants of Theorem \ref{theorem2d} to invariants.) However, there is one more twist, which is not seen in the cases of $X=\GL_n\backslash \PGL_{n+1}$, but is needed in the general case. Namely, instead of $\Gm$-invariants, one needs \emph{twisted} invariants with respect to a character of $\Gm$ (that is, of $F^\times$).

\subsubsection{}
To introduce this final piece of the puzzle, we recall from \cite{SaRankone} that the space $X\times X$ has two closed $G$-orbits of codimension larger than $1$: the diagonal $X^\diag$ (whose codimension we keep denoting by $d$), and a second closed $G$-orbit, whose codimension we will denote by $d'$. 
We define a character of $\Gm$ by $\chi_{d'}: z\mapsto |z|^{-\frac{d'}{2}+1}$. We will then understand the space of test densities for $(X\times X)/G$ as a subspace of the $(\Gm, \chi_{d'})$-equivariant elements of $D_\hor(\widetilde{J_X}^\circ, L_X)$, by composing \eqref{lifttotilde} with multiplication by 
\[\label{multfactor} |y|^{-\frac{d'}{2}+1} = |z(t_0 + \tau t_1)|^{-\frac{d'}{2}+1},\]
obtaining a map 
\begin{equation}\label{densityquotient-T}
 \mathcal D(X\times X) \to D_\hor(\widetilde{J_X}^\circ, L_X)^{(\Gm,\chi_{d'})}.
\end{equation}
The image $\mathcal D(\mathfrak X)$
is identified, as before, with the space $\mathcal S(\mathfrak X)^\st$ of stable test measures, if we fix an invariant measure on $X\times X$. To summarize, the map $\mathcal S(\mathfrak X)^\st \to \mathcal D(\mathfrak X)$ takes a measure $f(c)$ to 
\begin{equation}\label{meas2dens-T}
f(c) (dc)^{-\frac{1}{2}} |y|^{-\frac{d'}{2}+1} |t_1|^{-\frac{d}{2}+1} (d^\times z)^\frac{1}{2},
\end{equation}
in the trivialization \eqref{trivial-X}, where $y, c$ are given by \eqref{2dim}, \eqref{1dim}.

\begin{remark*}
The most convincing argument for the relevance of the character $\chi_{d'}$ is \cite[Proposition 6.1.5]{SaRankone}, describing orbital integrals in the neighborhood of $c=xy=0$ in terms of $\Gm$-orbital integrals on the $(x,y)$-plane, twisted by this character. However, a more conceptual understanding of it would be highly desirable.
\end{remark*}

\subsubsection{}
We also replace half-densities for the Kuznetsov formula of $G'=\PGL_2$ by half-densities for $\GL_2$ with central character $\chi_{d'}$, namely, we define an embedding
\[ \mathcal D(\mathfrak Y) \hookrightarrow D_\ver(\widetilde{J_X}^\circ, L_X)^{(\Gm,\chi_{d'})} \]
simply by multiplying the embedding $\mathcal D(\mathfrak Y) \hookrightarrow D_\ver(J_X^\circ, L_Y)$ of \S\ref{LYnatural} by the factor $|z|^{-\big(\frac{d'}{2}-1\big)}(d^\times z)^\frac{1}{2}$, and use the same  trivialization \eqref{trivial-Y} to identify $L_Y\simeq \underline{\mathbb C} \simeq L_X$.

The main result \cite[Theorem 1.3.1]{SaRankone} for the transfer operator in this case can now be formulated as follows:

\begin{theorem*}
There is an injective operator
 \[ \mathcal T: \mathcal D(\mathfrak Y) \to \mathcal D(\mathfrak X)\]
 given by integration along the leaves of the vertical foliation: 
 \[D_\ver (\widetilde{J_X}^\circ, L_X) \dashedrightarrow D_\hor (\widetilde{J_X}^\circ, L_X).
\]

Its inverse $\mathcal T^{-1}$, valued in an enlargement $\mathcal D_{L_X}^-(\mathfrak Y)\supset \mathcal D(\mathfrak Y)$ described in \cite[\S1.3]{SaRankone}, is given by integration along the leaves of the horizontal foliation:
\[D_\hor (\widetilde{J_X}^\circ, L_X) \dashedrightarrow D_\ver (\widetilde{J_X}^\circ, L_X),
\]
\end{theorem*}

\begin{proof}
 Let us compute the integral of a vertical half-density $\varphi \in D_\ver (\widetilde{J_X}^\circ, L_X)^{(\Gm,\chi_{d'})}$ along the vertical leaf passing through a point $j\in \widetilde{J_X}^\circ$ with coordinates $(t_0, t_1, \xi, z, \tau) \mod \mu_2$. If we think of it as a scalar-valued half-density, using our trivialization for $L_X$, the integral reads
\[ \underset{\nabla^X}{\underbrace{|t_1|^{-\frac{d}{2}+1}}} \int_{j'\in \mathscr F_{\ver, (x,y)}}  \underset{\nabla^X}{\underbrace{|t_1'|^{\frac{d}{2}-1}}} \cdot \varphi(j') \cdot |\omega \wedge \omega|^\frac{1}{2}.\]

As in the proof of Theorem \ref{mainthm-G}, we will use the 
Kostant section $(\xi,\tau)\mapsto \begin{pmatrix} \frac{\tau}{2} & \frac{\xi}{4} \\ 1 & \frac{\tau}{2} \end{pmatrix}$ to define a map
$\widetilde{J_X} \to T^*\tilde Y = \widetilde{\mathfrak c_X}^* \times \tilde G'$, 
where $\tilde G'=\GL_2$. One computes that the composition with projection to $\tilde G'$ is the map
\[ (t_0,t_1,\xi,z,\tau)\mapsto  \begin{pmatrix} z \\ & z \end{pmatrix} \begin{pmatrix} t_0 & \frac{\xi}{4}\cdot 2t_1 \\ 2t_1 & t_0\end{pmatrix} = \begin{pmatrix} 1 & \frac{t_0}{2t_1} \\ & 1 \end{pmatrix} \begin{pmatrix} & -z(2t_1)^{-1} \\ 2zt_1  \end{pmatrix}\begin{pmatrix} 1 & \frac{t_0}{2t_1} \\ & 1 \end{pmatrix} \]
Setting $a = zt_1$ and $b = z^{-1} t_1$, this can be rewritten as 
\[\begin{pmatrix} 1 & \frac{a^{-1}y+b^{-1}x}{4} \\ & 1 \end{pmatrix} \begin{pmatrix} & -(2b)^{-1} \\ 2a  \end{pmatrix}\begin{pmatrix} 1 & \frac{a^{-1}y+b^{-1}x}{4} \\ & 1 \end{pmatrix}.\]

Assume that $\varphi \in \mathcal D(\mathfrak Y)$ is the image of a measure $h = \Phi(t_1^2) d(t_1^2) \in \mathcal S(\mathfrak Y)$ (up to a choice of measure on $Y\times Y$). Then, with the trivialization \eqref{trivial-Y} composed with multiplication by $|z|^{-\big(\frac{d'}{2}-1\big)} (d^\times z)^\frac{1}{2}$, we have 
\begin{align*} \varphi(t_0,t_1,\xi,z,\tau) &= \psi\Big(\frac{t_0}{t_1} \Big) \frac{h}{|t_1 dt_1|^\frac{1}{2}}  |z|^{-\big(\frac{d'}{2}-1\big)} (d^\times z)^\frac{1}{2} 
\\ &= \psi\Big(\frac{t_0}{t_1} \Big)  \Phi(t_1^2) |t_1|  |z|^{-\big(\frac{d'}{2}-1\big)} d^\times(zt_1)^\frac{1}{2} d^\times(z^{-1}t_1)^\frac{1}{2} \\
&= \psi\Big(\frac{a^{-1}y+b^{-1}x}{2}\Big) \Phi(ab)  |ab|^\frac{1}{2} \Big|\frac{a}{b}\Big|^{-\frac{1}{2}\big(\frac{d'}{2}-1\big)} d^\times a^\frac{1}{2} d^\times b^\frac{1}{2}. 
\end{align*}

The half-density $|\omega \wedge \omega|^\frac{1}{2}$ can also be written $|4ab|^{-\frac{1}{2}} d^\times a^\frac{1}{2} d^\times b^\frac{1}{2} dx^\frac{1}{2} dy^\frac{1}{2}$; therefore, the integral above reads 
\begin{align*}
 |2|^{-1} dx^\frac{1}{2} dy^\frac{1}{2} |t_1|^{-\frac{d}{2}+1} \int  \Phi(ab)\cdot |ab|^{\frac{1}{2}\big(\frac{d}{2}-1\big)}\cdot \Big|\frac{a}{b}\Big|^{-\frac{1}{2}\big(\frac{d'}{2}-1\big)}  \cdot \psi\Big(\frac{a^{-1}y+b^{-1}x}{2}\Big) \, d^\times a \, d^\times b \\
 \approx dx^\frac{1}{2} dy^\frac{1}{2} |x|^{\frac{1}{4}(d+d')-1} |y|^{\frac{1}{4}(d-d')} |t_1|^{-\frac{d}{2}+1} \int \Phi\Big(\frac{xy}{4uv}\Big) |u|^{-\frac{1}{4}(d+d')+1} |v|^{-\frac{1}{4}(d-d')} \psi(u+v) \, d^\times u \, d^\times v,
\end{align*}
where $\approx$ means equality up to a constant that is independent of $\varphi$. By \eqref{meas2dens-T}, if this came from a measure in $\mathcal S(\mathfrak X)^\st$ (in the variable $c=xy$), that measure would be 
\[ dc \cdot |c|^{\frac{1}{4}(d+d')-1}  \int \Phi\Big(\frac{c}{4uv}\Big) |u|^{-\frac{1}{4}(d+d')+1} |v|^{-\frac{1}{4}(d-d')} \psi(u+v) \, d^\times u \, d^\times v.\]
This is the formula \cite[(9)]{SaRankone} for the transfer operator. (Note that $\frac{1}{4}(d+d') -1 = s_1- \frac{1}{2}$ and $\frac{1}{4}(d-d') = s_2- \frac{1}{2}$ in the notation of \cite[(83), (85)]{SaRankone}, and that the coordinates used in that paper are $4t_1^2$ to describe measures in $\mathcal S(\mathfrak Y)$, and $c$ to describe measures in $\mathcal S(\mathfrak X)^\st$.)

Again, I leave it to the reader to check the formula for $\mathcal T^{-1}$.
\end{proof}

\subsubsection{} \label{nonsplit}
I finish this section with a brief discussion of a case where a section $\mathfrak c_X^*\to T^*X$ does not exist. Let $X=T\backslash \PGL_2$, where $T$ is a nonsplit torus, splitting over a quadratic extension $E/F$. In this case, the transfer operator 
\[ \mathcal T: \mathcal S(\mathfrak Y) \to \mathcal S(\mathfrak X)\]
was computed in \cite{SaBE1}, and can be described as follows: 

Instead of defining $\widetilde{J_X}$ to be the group scheme of regular centralizers in $\GL_2$, define it to be the $\operatorname{Gal}(E/F)$-twist of that, determined by the automorphism $(t_0,t_1,\xi,z,\tau)\mapsto (t_0, t_1, \xi, z^{-1}, -\tau)$; that is, $\widetilde{J_X}$ will be isomorphic to the group scheme of regular centralizers in the quasisplit unitary group $U_2$. The transfer operator, now, it obtained as in Theorem \ref{theorem2d}, by descending the operator of integration along the leaves of the corresponding ``vertical'' foliation on $\widetilde{J_X}$ to $U_1$-coinvariants.

Both the Schwartz space and the descent to $U_1$-coinvariants, here, need to be understood in a sophisticated, ``stacky'' way. 
Namely, the full space $\mathcal S(\mathfrak X)$ includes a ``pure inner form'' as in \eqref{Schwartz},
\[ \mathcal S(\mathfrak X) = \mathcal S(X\times X)_G \oplus \mathcal S(X^\alpha\times X^\alpha)_{G^\alpha}, \]
where $X^\alpha \simeq T\backslash G^\alpha$, with $G^\alpha =PD^\times$, the projective multiplicative group of the quaternion division algebra. 
Similarly, the $U_1$-coinvariants of $D_\ver (\widetilde{J_X}^\circ, L_X)$, $D_\hor (\widetilde{J_X}^\circ, L_X)$ need to be understood in a stacky way. Explicitly, recall that in the split case the space $D_\hor (\widetilde{J_X}^\circ, L_X)$ was the space of half-densities on the $(x,y)$-plane (in coordinates \eqref{2dim}), with $\Gm$ acting as $z\cdot (x,y) = (z^{-1}x, zy)$. In the nonsplit case, the $(x,y)$-plane becomes the space 
$V=\operatorname{Res}_{E/F} \Ga$, and instead of $U_1$-coinvariants of the space $\mathcal D(V(F))$, one needs to consider the direct sum 
\[ \mathcal D(V(F))_{U_1} \oplus \mathcal D(V'(F))_{U_1},\]
where $V'$ is the twist of $V$ by the nontrivial $U_1$-torsor. The same interpretation is needed for ``stacky'' $U_1$-coinvariants of the space $D_\ver (\widetilde{J_X}^\circ, L_X)$ (with the coordinates $a = zt_1$, $b=z^{-1}t_1$ for the leaves of the horizontal foliation interchanged by the Galois action), and the operator of ``integration along vertical half-densities'' -- essentially, a Fourier transform from the Galois-twisted $(x,y)$-plane to the Galois-twisted $(a^{-1}, b^{-1})$-plane -- naturally descends to give the transfer operator $\mathcal T$ of \cite[Theorem 5.1]{SaBE1}.

\section{Problems for the near future} \label{sec:future}

\subsection{The relative Langlands conjectures}

The relative Langlands conjectures presented in \S\ref{conj-local}--\ref{conj-global} do not have the precision of the local conjectures of Gan--Gross--Prasad \cite{GGP}, or the global conjectures of Ichino--Ikeda \cite{II}. Moreover, an extension of those conjectures to Arthur packets not appearing in the $L^2$-decomposition is required, as in \cite{GGP-nontempered}. 

It is therefore an important problem to refine the existing conjectures. It is also a fascinating one: as always, finding a way to blend several known cases into a uniform theory can lead to new insights about the nature of the problems. The geometric relative Langlands conjectures proposed in joint work with Ben-Zvi and Venkatesh can probably assist in this direction, providing a geometric spectral answer to automorphic problems, which can then be translated to number theory by de-categorifying.

\subsection{Transfer operators in higher rank}

The most important problem ``beyond endoscopy,'' in my view, is to understand transfer operators in higher rank, and for morphisms of $L$-groups ${^LG_X}\to {^LG_Y}$ that are not isomorphisms. Regarding the latter, despite the traditional emphasis on the Arthur--Selberg trace formula, it might be better, as first observed by Sarnak \cite{Sarnak}, to try to compare Kuznetsov formulas, which, according to the results of \cite{SaRankone}, seem to be the ``base cases'' for every comparison. This can be ``explained'' by the simple structure \eqref{fromJX-Kuz} of the Kuznetsov cotangent stack.

If the ideas discussed in this paper have any merit, understanding transfer operators in terms of ``quantization'' would involve several steps, including the following:

\subsubsection{} Develop a theory of ``geometric quantization'' for derived symplectic stacks, whose output includes the Schwartz spaces of stacks defined in \cite{SaStacks}. The cases presented here, and in particular the construction of the maps \eqref{densityquotient}, \eqref{densityquotient-T}, could provide some hints on how to do that, but the various twists involved need to be better understood. 

\subsubsection{} \label{cotangentextension} Obtain a better understanding of the structure of coisotropic Hamiltonian spaces and the cotangent stacks appearing in the relative trace formula. The diagram \eqref{fromJX}, arising from the work of Knop, which was interpreted as a ``birational'' description of  $T^*\mathfrak X$, does not capture the difference between $T^*\mathfrak X$, and the Kuznetsov cotangent stack with the same $L$-group. This difference seems to be significant for the structure of transfer operators and for spaces of test measures. For example, the enlarged spaces $\mathcal S_{L_X}^-(\mathfrak Y)$ of test measures for the Kuznetsov formula in Theorems \ref{mainthm-G}, \ref{mainthm-T} should be seen as quantizations of $T^*\mathfrak X$, which is strictly larger than the Kuznetsov stack $T^*\mathfrak Y$. This difference can also explain ``Galois twists'' of transfer operators, as in the example of \S\ref{nonsplit}, where $T^*\mathfrak X$ failed to admit an analog of a Kostant section.

A baby case of the idea that embeddings of Hamiltonian spaces correspond to enlarged Schwartz spaces is Iwasawa--Tate theory, where the embedding $T^*\Gm\hookrightarrow T^*\mathbb A^1$ corresponds to the enlargement $\mathcal S(F^\times)\hookrightarrow \mathcal S(F)$. This can be generalized to toric stacks, as described in \cite[\S5.2]{Ngo-Takagi}. 
For a torus $A$, a collection of coweights $\mu:\Gm^r \to A$ defines a stack $\overline{A}_\mu:= \mathbb A^r/\ker(\mu)$ with an action of $A$. The \emph{dual Hamiltonian space} of such a stack is defined as the symplectic $\check A$-vector space $\check M$ with weights $\{\pm \mu\}$. 

To give an example of the descriptions of Hamiltonian spaces envisioned here, in ongoing joint work with Ben-Zvi and Venkatesh we take a step beyond Knop's theory, modeling the most regular locus of $M=T^*X$, up to codimension $2$, on a the analog of a ``toric stack'' for the group scheme $J_X$. Our result confirms an observation of V.\ Lafforgue, shared in private communication several years ago. For example, for spherical varieties $X$ with $\check G_X=\check G$ and $\check A_X=\check A$, our description uses the toric stack $\overline{A_X}$ corresponding to the dual Hamiltonian space $\check M=V_X$ of the spherical variety, described in terms of its colors in \S\ref{density-smooth}.

\begin{theorem*} \label{Hamiltonian}
In the setting above, there is an action of $W_X$ on the toric cotangent stack $T^*{\overline{A_X}}$, and an open dense subset ${\mathfrak c_X^*}'\subset \mathfrak c_X^*$, whose complement has codimension $\ge 2$, such that the restriction 
$M'\subset M$ to the image of ${\mathfrak c_X^*}'$ under a Kostant section $\kappa: \mathfrak c_X^*\to \mathfrak g^*$ admits a $J_X$-equivariant symplectomorphism
 \[ M' \simeq  \Big( \operatorname{Res}_{\mathfrak a_X^*/\mathfrak c_X^*} T^*\overline{A_X}  \Big)^{W_X} \times_{\mathfrak c_X^*} {\mathfrak c_X^*}'.\]
\end{theorem*}

\subsubsection{} Use cotangent spaces to understand transfer operators. As we saw in the discussion of cases of type $T$ in \S\ref{sstypeT}, a ``naive'' change of geometric quantization on $J_X$ did not give the correct transfer operators; instead, one has to pass to the group scheme $\widetilde{J_X}$ associated with the root datum of $\GL_2$, producing a $2$-dimensional Fourier transform. 

Ongoing joint work with C.\ Wan, comparing the Kuznetsov to the Arthur--Selberg trace formula for $\GL_n$, suggests that, in higher-rank cases, the transfer operator for a comparison to a Kuznetsov quotient with the same dual group might be given by an $r$-dimensional integration, where $r$ is, roughly, half the dimension of the non-zero weight spaces of the representation $V_X$ (\S\ref{Lconjecture}) of the dual group.  In view of the discussion of \S\ref{cotangentextension}, this seems to be closely related to the structure of $T^*\mathfrak X$. For example, in the setting of Theorem \ref{cotangentextension} (such as in the case of Gan--Gross--Prasad periods), one could speculate that the transfer operator to the Kuznetsov formula is somehow a ``descent'' of a Fourier transform in $\frac{\dim \check M}{2}$-dimensions. In this ``dream,'' the following three objects would be closely related: the $L$-value $L_X$  associated to a spherical variety (encoded in a dual Hamiltonian space $\check M$), the fine structure of the Hamiltonian space $M=T^*X$ (whose quantization is the space of test measures on $X$), and the transfer operator between the relative trace formula of $X$ and the Kuznetsov formula with the same dual group.

\subsection{Poisson summation formula} 

Understanding the local transfer operators should be followed by a global comparison of trace formulas. For example, comparing the relative trace formula for any spherical variety $X$ with the Kuznetsov formula with the same dual group should amount to commutativity of the diagram
\begin{equation}\label{RTFdiagram} 
\xymatrix{ \mathcal S_{L_X}^-(\mathfrak Y(\mathbb A)) \ar[rr]^{\mathcal T}\ar[dr]_{\RTF_{\mathfrak Y}} & & \ar[dl]^{\RTF_{\mathfrak X}}\mathcal S(\mathfrak X(\mathbb A))^\st\\
& \mathbb C &},
\end{equation}
where $\mathcal S_{L_X}^-(\mathfrak Y(\mathbb A))$ is a suitable enlarged space of test measures for the Kuznetsov formula, related to the $L$-value $L_X$ of $X$. 

Having a formula for the transfer operator in terms of Fourier transforms (as in Theorems \ref{mainthm-G}, \ref{mainthm-T}) gives hope of employing the Poisson summation formula to establish commutativity of \eqref{RTFdiagram}. However, this is far from straightforward, as the spaces of stable test measures are nonstandard. In Altu\u{g}'s work \cite{Altug1}, the approximate functional equation was used for the trace formula of $X=\GL_2$, obtaining an expression similar to the Kuznetsov formula (in particular, containing Kloosterman sums), but not quite equal to it. 

A different approach was introduced in \cite{SaBE2}, for the case $X=T\backslash \PGL_2$. It is based on the idea of \emph{deforming} spaces of test measures and transfer operators in analytic families depending on a parameter $s$ (which moves the point of evaluation of $L_X$), so that in some domain for $s$ the Poisson summation formula is valid. It is likely that this method can be applied more generally, but it requires a better understanding of the idea of deforming spaces of test measures (orbital integrals).

\subsection{Hankel transforms} \label{Hankel}

In the recent literature on automorphic forms, the term ``Hankel transforms'' has been used to describe two distinct conjectural notions:
\begin{itemize}[leftmargin=2em]
\item The nonlinear Fourier transforms $\mathcal S_\rho(G(F)) \to \mathcal S_{\rho^*}(G(F))$ between nonstandard spaces of Schwartz functions (or measures) on a reductive group over a local field, adapted to a representation $\rho$ of its dual group, and its dual. These spaces and operators would generalize Fourier transform on the space $\mathcal S(\operatorname{Mat}_n(F))$ of Godement--Jacquet theory, for the case $\rho=$ the standard representation of $\check G=\GL_n$, and would similarly give rise to the local functional equation for the $L$-functions associated to $\rho$. They were introduced by Braverman and Kazhdan \cite{BKgamma1,BKgamma2}, and advanced in the work of Ng\^o \cite{Ngo-Takagi}. 

\item The descent of such transforms to spaces of test measures for the Arthur--Selberg trace formula, or for the Kuznetsov formula. In the latter case, those would be operators 
\[ \mathcal H_\rho: \mathcal S_\rho(\mathfrak Y) \to \mathcal S_{\rho^*} (\mathfrak Y)\] 
between enlarged spaces of measures for the Kuznetsov formula, such as those encountered in Theorems \ref{mainthm-G}, \ref{mainthm-T}. 
\end{itemize}

The two notions are closely related, but it is the latter that we would like to focus on here. It is natural to ask the question of whether one can describe $\mathcal H_\rho$ explicitly, and prove a Poisson summation formula, globally, in the sense that the diagram 
\begin{equation}\label{Hankeldiagram} 
\xymatrix{ \mathcal S_\rho(\mathfrak Y(\mathbb A)) \ar[rr]^{\mathcal H_\rho}\ar[dr]_{\RTF_{\mathfrak Y}} & & \ar[dl]^{\RTF_{\mathfrak Y}}\mathcal S_{\rho^*}(\mathfrak Y(\mathbb A))\\
& \mathbb C &}
\end{equation}
should commute. This would lead to an independent proof of the functional equation of the pertinent $L$-functions. 

Such Hankel transforms have been described by Jacquet \cite{Jacquet-smoothtransfer} for $\rho=$ the standard representation of $\GL_n$ (the paper \cite{Herman} is closely related), and by me \cite{SaTransfer2} for $\rho=$ the symmetric square representation of $\GL_2$.  It would be interesting to examine if these formulas admit an interpretation in terms of quantization, like the transfer operators in this paper.

It seems counter to the strategy of the Langlands program to seek such a proof of the functional equation, independent of functoriality. On the other hand, the similarity between diagrams \eqref{RTFdiagram} and \eqref{Hankeldiagram} is enticing. More fundamentally, trace formulas with nonstandard test functions, such as those in Langlands' original ``beyond endoscopy'' proposal, or the Kuznetsov formula appearing in \eqref{RTFdiagram}, have the $L$-functions embedded into them. Obtaining the spectral decomposition of those formulas will likely require more than ``brute force'' analytic number theory, and a Poisson summation formula of the form \eqref{Hankeldiagram} could help resolve the problem. This idea was successfully employed in \cite{SaBE2} for a new proof of Waldspurger's formula for toric periods in $\PGL_2$ via a nonstandard comparison of the form \eqref{RTFdiagram}. Therefore, it might be that, in the ``beyond endoscopy'' program, functoriality and the functional equation of $L$-functions should be studied hand-in-hand.


\section*{acknowledgments}
My understanding of the subject has been greatly influenced by  conversations with David Ben-Zvi, Rapha\"el Beuzart-Plessis, Ng\^o Bao Ch\^au, Akshay Venkatesh, Chen Wan, and Jonathan Wang. I thank them all for generously sharing their ideas in our exciting mathematical journeys.

This work was partially supported by NSF grants DMS-1939672 and DMS-2101700, and by a Simons Fellowship in Mathematics.


\bibliographystyle{emss}
\bibliography{sakellaridis}

\begin{thebibliography}{10}
\providecommand{\url}[1]{\texttt{#1}}
\providecommand{\urlprefix}{URL }
\providecommand{\eprint}[2][]{\url{#2}}

\bibitem{AGSchwartz}
A.~Aizenbud and D.~Gourevitch, Schwartz functions on {N}ash manifolds.
  \emph{Int. Math. Res. Not. IMRN}  (2008), no.~5, Art. ID rnm 155, 37.
  \MR{2418286}

\bibitem{Altug1}
S.~A. Altu\u{g}, Beyond endoscopy via the trace formula: 1. {P}oisson summation
  and isolation of special representations. \emph{Compos. Math.} \textbf{151}
  (2015), no.~10, 1791--1820. \MR{3414386}

\bibitem{Altug2}
S.~A. Altu\u{g}, Beyond endoscopy via the trace formula, {II}: {A}symptotic
  expansions of {F}ourier transforms and bounds towards the {R}amanujan
  conjecture. \emph{Amer. J. Math.} \textbf{139} (2017), no.~4, 863--913.
  \MR{3689319}

\bibitem{Altug3}
S.~A. Altu\u{g}, Beyond endoscopy via the trace formula---{III} {T}he standard
  representation. \emph{J. Inst. Math. Jussieu} \textbf{19} (2020), no.~4,
  1349--1387. \MR{4120811}

\bibitem{AYY}
J.~An, J.-K. Yu, and J.~Yu, On the dimension datum of a subgroup and its
  application to isospectral manifolds. \emph{J. Differential Geom.}
  \textbf{94} (2013), no.~1, 59--85. \MR{3031860}

\bibitem{Arthur-invariant}
J.~Arthur, The trace formula in invariant form. \emph{Ann. of Math. (2)}
  \textbf{114} (1981), no.~1, 1--74. \MR{625344}

\bibitem{Beuzart-GP}
R.~Beuzart-Plessis, A local trace formula for the {G}an-{G}ross-{P}rasad
  conjecture for unitary groups: the {A}rchimedean case. \emph{Ast\'{e}risque}
  (2020), no. 418, viii + 299. \MR{4146145}

\bibitem{Beuzart-unitary}
R.~Beuzart-Plessis, Multiplicities and {P}lancherel formula for the space of
  nondegenerate {H}ermitian matrices. \emph{J. Number Theory} \textbf{230}
  (2022), 5--63. \MR{4327948}

\bibitem{BCZ}
R.~Beuzart-Plessis, P.-H. Chaudouard, and M.~Zydor, The global
  {G}an-{G}ross-{P}rasad conjecture for unitary groups: the endoscopic case.
  \emph{Publ. Math. Inst. Hautes \'{E}tudes Sci.} \textbf{135} (2022),
  183--336. \MR{4426741}

\bibitem{BLZZ}
R.~Beuzart-Plessis, Y.~Liu, W.~Zhang, and X.~Zhu, Isolation of cuspidal
  spectrum, with application to the {G}an-{G}ross-{P}rasad conjecture.
  \emph{Ann. of Math. (2)} \textbf{194} (2021), no.~2, 519--584. \MR{4298750}

\bibitem{BezFin}
R.~Bezrukavnikov and M.~Finkelberg, Equivariant {S}atake category and
  {K}ostant-{W}hittaker reduction. \emph{Mosc. Math. J.} \textbf{8} (2008),
  no.~1, 39--72, 183. \MR{2422266}

\bibitem{BNS}
A.~Bouthier, B.~C. Ng\^o, and Y.~Sakellaridis, On the formal arc space of a
  reductive monoid. \emph{Amer. J. Math.} \textbf{138} (2016), no.~1, 81--108.

\bibitem{BFGM}
A.~Braverman, M.~Finkelberg, D.~Gaitsgory, and I.~Mirkovi\'{c}, Intersection
  cohomology of {D}rinfeld's compactifications. \emph{Selecta Math. (N.S.)}
  \textbf{8} (2002), no.~3, 381--418. \MR{1931170}

\bibitem{BK1}
A.~Braverman and D.~Kazhdan, On the {S}chwartz space of the basic affine space.
  \emph{Selecta Math. (N.S.)} \textbf{5} (1999), no.~1, 1--28. \MR{1694894}

\bibitem{BKgamma1}
A.~Braverman and D.~Kazhdan, {$\gamma$}-functions of representations and
  lifting. pp. 237--278, Special Volume, Part I, 2000. \MR{1826255}

\bibitem{BK2}
A.~Braverman and D.~Kazhdan, Normalized intertwining operators and nilpotent
  elements in the {L}anglands dual group. pp. 533--553, 2, 2002. \MR{1988971}

\bibitem{BKgamma2}
A.~Braverman and D.~Kazhdan, {$\gamma$}-sheaves on reductive groups. In
  \emph{Studies in memory of {I}ssai {S}chur ({C}hevaleret/{R}ehovot, 2000)},
  pp. 27--47, Progr. Math. 210, Birkh\"{a}user Boston, Boston, MA, 2003.
  \MR{1985192}

\bibitem{Brion}
M.~Brion, Vers une g\'en\'eralisation des espaces sym\'etriques. \emph{J.
  Algebra} \textbf{134} (1990), no.~1, 115--143. \MR{1068418}

\bibitem{Bump-RS}
D.~Bump, The {R}ankin-{S}elberg method: an introduction and survey. In
  \emph{Automorphic representations, {$L$}-functions and applications: progress
  and prospects}, pp. 41--73, Ohio State Univ. Math. Res. Inst. Publ. 11, de
  Gruyter, Berlin, 2005. \MR{2192819}

\bibitem{Casselman}
W.~Casselman, The unramified principal series of {$p$}-adic groups. {I}. {T}he
  spherical function. \emph{Compositio Math.} \textbf{40} (1980), no.~3,
  387--406. \MR{571057}

\bibitem{CS}
W.~Casselman and J.~Shalika, The unramified principal series of {$p$}-adic
  groups. {II}. {T}he {W}hittaker function. \emph{Compositio Math.} \textbf{41}
  (1980), no.~2, 207--231. \MR{581582}

\bibitem{Delorme-real}
P.~Delorme, Formule de {P}lancherel pour les espaces sym\'etriques r\'eductifs.
  \emph{Ann. of Math. (2)} \textbf{147} (1998), no.~2, 417--452. \MR{1626757
  (99d:22022)}

\bibitem{Delorme-Plancherel}
P.~Delorme, Neighborhoods at infinity and the {P}lancherel formula for a
  reductive {$p$}-adic symmetric space. \emph{Math. Ann.} \textbf{370} (2018),
  no. 3-4, 1177--1229. \MR{3770165}

\bibitem{FLO}
B.~Feigon, E.~Lapid, and O.~Offen, On representations distinguished by unitary
  groups. \emph{Publ. Math. Inst. Hautes \'{E}tudes Sci.} \textbf{115} (2012),
  185--323. \MR{2930996}

\bibitem{FLN}
E.~Frenkel, R.~Langlands, and B.~C. Ng\^{o}, Formule des traces et
  fonctorialit\'{e}: le d\'{e}but d'un programme. \emph{Ann. Sci. Math.
  Qu\'{e}bec} \textbf{34} (2010), no.~2, 199--243. \MR{2779866}

\bibitem{GN}
D.~Gaitsgory and D.~Nadler, Spherical varieties and {L}anglands duality.
  \emph{Mosc. Math. J.} \textbf{10} (2010), no.~1, 65--137, 271. \MR{2668830}

\bibitem{GGP}
W.~T. Gan, B.~H. Gross, and D.~Prasad, Symplectic local root numbers, central
  critical {$L$} values, and restriction problems in the representation theory
  of classical groups. pp. 1--109, 346, 2012. \MR{3202556}

\bibitem{GGP-nontempered}
W.~T. Gan, B.~H. Gross, and D.~Prasad, Branching laws for classical groups: the
  non-tempered case. \emph{Compos. Math.} \textbf{156} (2020), no.~11,
  2298--2367. \MR{4190046}

\bibitem{Gan-Ichino}
W.~T. Gan and A.~Ichino, The {G}ross-{P}rasad conjecture and local theta
  correspondence. \emph{Invent. Math.} \textbf{206} (2016), no.~3, 705--799.
  \MR{3573972}

\bibitem{GQT}
W.~T. Gan, Y.~Qiu, and S.~Takeda, The regularized {S}iegel-{W}eil formula (the
  second term identity) and the {R}allis inner product formula. \emph{Invent.
  Math.} \textbf{198} (2014), no.~3, 739--831. \MR{3279536}

\bibitem{Gan-Wan}
W.~T. Gan and X.~Wan, Relative character identities and theta correspondence.
  In \emph{Relative trace formulas}, edited by W.~M{\"u}ller, S.~W. Shin, and
  N.~Templier, pp. 101--186, Springer International Publishing, Cham, 2021.

\bibitem{GHL}
J.~Getz, C.-H. Hsu, and S.~Leslie, Harmonic analysis on certain spherical
  varieties, preprint, \eprint{2103.10261}.

\bibitem{GL1}
J.~R. Getz and B.~Liu, A summation formula for triples of quadratic spaces.
  \emph{Adv. Math.} \textbf{347} (2019), 150--191. \MR{3916514}

\bibitem{GL2}
J.~R. Getz and B.~Liu, A refined {P}oisson summation formula for certain
  {B}raverman-{K}azhdan spaces. \emph{Sci. China Math.} \textbf{64} (2021),
  no.~6, 1127--1156. \MR{4268887}

\bibitem{HC3}
Harish-Chandra, Harmonic analysis on real reductive groups. {III}. {T}he
  {M}aass-{S}elberg relations and the {P}lancherel formula. \emph{Ann. of Math.
  (2)} \textbf{104} (1976), no.~1, 117--201. \MR{0439994 (55 \#12875)}

\bibitem{Hecke1}
E.~Hecke, \"{U}ber die {B}estimmung {D}irichletscher {R}eihen durch ihre
  {F}unktionalgleichung. \emph{Math. Ann.} \textbf{112} (1936), no.~1,
  664--699. \MR{1513069}

\bibitem{Hecke2}
E.~Hecke, \"{U}ber {M}odulfunktionen und die {D}irichletschen {R}eihen mit
  {E}ulerscher {P}roduktentwicklung. {I}. \emph{Math. Ann.} \textbf{114}
  (1937), no.~1, 1--28. \MR{1513122}

\bibitem{Hecke3}
E.~Hecke, \"{U}ber {M}odulfunktionen und die {D}irichletschen {R}eihen mit
  {E}ulerscher {P}roduktentwicklung. {II}. \emph{Math. Ann.} \textbf{114}
  (1937), no.~1, 316--351. \MR{1513142}

\bibitem{Herman}
P.~E. Herman, The functional equation and beyond endoscopy. \emph{Pacific J.
  Math.} \textbf{260} (2012), no.~2, 497--513. \MR{3001802}

\bibitem{II}
A.~Ichino and T.~Ikeda, On the periods of automorphic forms on special
  orthogonal groups and the {G}ross-{P}rasad conjecture. \emph{Geom. Funct.
  Anal.} \textbf{19} (2010), no.~5, 1378--1425. \MR{2585578}

\bibitem{Iwasawa}
K.~Iwasawa, A note on functions. In \emph{Proceedings of the {I}nternational
  {C}ongress of {M}athematicians, {C}ambridge, {M}ass., 1950, vol. 1}, p. 322,
  Amer. Math. Soc., Providence, R. I., 1952.

\bibitem{Jacquet-symmetric}
H.~Jacquet, Automorphic spectrum of symmetric spaces. In \emph{Representation
  theory and automorphic forms ({E}dinburgh, 1996)}, pp. 443--455, Proc.
  Sympos. Pure Math. 61, Amer. Math. Soc., Providence, RI, 1997. \MR{1476509}

\bibitem{Jacquet-factorization}
H.~Jacquet, Factorization of period integrals. \emph{J. Number Theory}
  \textbf{87} (2001), no.~1, 109--143. \MR{1816039}

\bibitem{Jacquet-smoothtransfer}
H.~Jacquet, Smooth transfer of {K}loosterman integrals. \emph{Duke Math. J.}
  \textbf{120} (2003), no.~1, 121--152. \MR{2010736}

\bibitem{JL}
H.~Jacquet and K.~F. Lai, A relative trace formula. \emph{Compositio Math.}
  \textbf{54} (1985), no.~2, 243--310. \MR{783512}

\bibitem{JLR}
H.~Jacquet, K.~F. Lai, and S.~Rallis, A trace formula for symmetric spaces.
  \emph{Duke Math. J.} \textbf{70} (1993), no.~2, 305--372. \MR{1219816}

\bibitem{Johnstone}
D.~Johnstone, A {G}elfand--{G}raev formula and stable transfer factors for
  $\operatorname{SL}_n(f)$, preprint, \eprint{1611.06291}.

\bibitem{Kirillov}
A.~A. Kirillov, \emph{Lectures on the orbit method}. Graduate Studies in
  Mathematics 64, American Mathematical Society, Providence, RI, 2004.
  \MR{2069175}

\bibitem{KnWeyl}
F.~Knop, Weylgruppe und {M}omentabbildung. \emph{Invent. Math.} \textbf{99}
  (1990), no.~1, 1--23. \MR{1029388}

\bibitem{KnAut}
F.~Knop, Automorphisms, root systems, and compactifications of homogeneous
  varieties. \emph{J. Amer. Math. Soc.} \textbf{9} (1996), no.~1, 153--174.
  \MR{1311823}

\bibitem{KnSch}
F.~Knop and B.~Schalke, The dual group of a spherical variety. \emph{Trans.
  Moscow Math. Soc.} \textbf{78} (2017), 187--216. \MR{3738085}

\bibitem{Kostant}
B.~Kostant, Quantization and unitary representations. In \emph{Lectures in
  modern analysis and applications, {III}}, pp. 87--208. Lecture Notes in
  Math., Vol. 170, 1970. \MR{0294568}

\bibitem{LLafforgue}
L.~Lafforgue, Noyaux du transfert automorphe de {L}anglands et formules de
  {P}oisson non lin\'{e}aires. \emph{Jpn. J. Math.} \textbf{9} (2014), no.~1,
  1--68. \MR{3173438}

\bibitem{Langlands-Euler}
R.~P. Langlands, \emph{Euler products}. Yale Mathematical Monographs 1, Yale
  University Press, New Haven, Conn.-London, 1971. \MR{0419366}

\bibitem{Langlands-BE}
R.~P. Langlands, Beyond endoscopy. In \emph{Contributions to automorphic forms,
  geometry, and number theory}, pp. 611--697, Johns Hopkins Univ. Press,
  Baltimore, MD, 2004. \MR{2058622 (2005f:11102)}

\bibitem{Langlands-ST}
R.~P. Langlands, Singularit\'es et transfert. \emph{Ann. Sci. Math. Qu\'{e}bec}
  \textbf{37} (2013), no.~2, 173--253. \MR{3117742}

\bibitem{LM-unitary}
E.~Lapid and Z.~Mao, On {W}hittaker-{F}ourier coefficients of automorphic forms
  on unitary groups: reduction to a local identity. In \emph{Advances in the
  theory of automorphic forms and their {$L$}-functions}, pp. 295--320,
  Contemp. Math. 664, Amer. Math. Soc., Providence, RI, 2016. \MR{3502987}

\bibitem{LM}
E.~Lapid and Z.~Mao, On an analogue of the {I}chino-{I}keda conjecture for
  {W}hittaker coefficients on the metaplectic group. \emph{Algebra Number
  Theory} \textbf{11} (2017), no.~3, 713--765. \MR{3649366}

\bibitem{WWLi}
W.-W. Li, The {W}eil representation and its character, 2008, {M}aster's thesis,
  Universiteit Leiden.

\bibitem{Maass}
H.~Maass, \"{U}ber eine neue {A}rt von nichtanalytischen automorphen
  {F}unktionen und die {B}estimmung {D}irichletscher {R}eihen durch
  {F}unktionalgleichungen. \emph{Math. Ann.} \textbf{121} (1949), 141--183.
  \MR{31519}

\bibitem{MR}
C.~M{\oe}glin and D.~Renard, S{\'e}ries discr{\`e}tes des espaces
  sym{\'e}triques et paquets d'{A}rthur, preprint, \eprint{1906.00725}.

\bibitem{MW-GP}
C.~M{\oe}glin and J.-L. Waldspurger, La conjecture locale de {G}ross-{P}rasad
  pour les groupes sp\'{e}ciaux orthogonaux: le cas g\'{e}n\'{e}ral. pp.
  167--216, 347, 2012. \MR{3155346}

\bibitem{Morimoto}
K.~Morimoto, On a certain local identity for {L}apid--{M}ao's conjecture and
  formal degree conjecture: even unitary group case, preprint,
  \eprint{1902.04910}.

\bibitem{Na-real}
D.~Nadler, Perverse sheaves on real loop {G}rassmannians. \emph{Invent. Math.}
  \textbf{159} (2005), no.~1, 1--73. \MR{2142332}

\bibitem{Ngo-FL}
B.~C. Ng\^o, Le lemme fondamental pour les alg\`ebres de {L}ie. \emph{Publ.
  Math. Inst. Hautes \'Etudes Sci.}  (2010), no. 111, 1--169. \MR{2653248}

\bibitem{Ngo-PS}
B.~C. Ng\^o, On a certain sum of automorphic {$L$}-functions. In
  \emph{Automorphic forms and related geometry: assessing the legacy of {I}.
  {I}. {P}iatetski-{S}hapiro}, pp. 337--343, Contemp. Math. 614, Amer. Math.
  Soc., Providence, RI, 2014. \MR{3220933}

\bibitem{Ngo-Takagi}
B.~C. Ng\^{o}, Hankel transform, {L}anglands functoriality and functional
  equation of automorphic {$L$}-functions. \emph{Jpn. J. Math.} \textbf{15}
  (2020), no.~1, 121--167. \MR{4068833}

\bibitem{Prasad-relative}
D.~Prasad, A `relative' local langlands correspondence, preprint,
  \eprint{1512.04347}.

\bibitem{Rankin}
R.~A. Rankin, Contributions to the theory of {R}amanujan's function $\tau(n)$
  and similar arithmetical functions: {II}. {T}he order of the {F}ourier
  coefficients of integral modular forms. \emph{Mathematical Proceedings of the
  Cambridge Philosophical Society} \textbf{35} (1939), no.~3, 357–372.

\bibitem{Riemann}
B.~Riemann, Ueber die {A}nzahl der {P}rimzahlen unter einer gegebenen
  {G}r\"osse. \emph{Monatsberichte der Berliner Akademie}  (November 1859).

\bibitem{Rudnick}
Z.~Rudnick, \emph{Poincare series}. Ph.D. thesis, 1990, {Y}ale University.
  \MR{2638802}

\bibitem{SaRS}
Y.~Sakellaridis, Spherical varieties and integral representations of
  {$L$}-functions. \emph{Algebra \& Number Theory} \textbf{6} (2012), no.~4,
  611--667.

\bibitem{SaBE1}
Y.~Sakellaridis, Beyond endoscopy for the relative trace formula {I}: local
  theory. In \emph{Automorphic representations and {L}-functions}, pp.
  521--590, Amer. Math. Soc., Providence, RI, 2013.

\bibitem{SaSph}
Y.~Sakellaridis, Spherical functions on spherical varieties. \emph{Amer. J.
  Math.} \textbf{135} (2013), no.~5, 1291--1381.

\bibitem{SaStacks}
Y.~Sakellaridis, The {S}chwartz space of a smooth semi-algebraic stack.
  \emph{Selecta Math. (N.S.)} \textbf{22} (2016), no.~4, 2401--2490.

\bibitem{SaHowe}
Y.~Sakellaridis, Plancherel decomposition of {H}owe duality and {E}uler
  factorization of automorphic functionals. In \emph{Representation theory,
  number theory, and invariant theory}, pp. 545--585, Progr. Math. 323,
  Birkh\"auser/Springer, Cham, 2017. \MR{3753923}

\bibitem{SaBE2}
Y.~Sakellaridis, Beyond endoscopy for the relative trace formula {II}: global
  theory. \emph{J. Inst. Math. Jussieu} \textbf{18} (2019), no.~2, 347--447.
  \MR{3915291}

\bibitem{SaRankone}
Y.~Sakellaridis, Functorial transfer between relative trace formulas in rank
  {$1$}. \emph{Duke Math. J.} \textbf{170} (2021), no.~2, 279--364.
  \MR{4202495}

\bibitem{SaTransfer1}
Y.~Sakellaridis, Transfer operators and {H}ankel transforms between relative
  trace formulas, {I}: {C}haracter theory. \emph{Adv. Math.} \textbf{394}
  (2022), Paper No. 108010, 75. \MR{4355722}

\bibitem{SaTransfer2}
Y.~Sakellaridis, Transfer operators and {H}ankel transforms between relative
  trace formulas, {II}: {R}ankin-{S}elberg theory. \emph{Adv. Math.}
  \textbf{394} (2022), Paper No. 108039, 104. \MR{4355733}

\bibitem{SV}
Y.~Sakellaridis and A.~Venkatesh, Periods and harmonic analysis on spherical
  varieties. \emph{Ast\'erisque}  (2017), no. 396, 360.

\bibitem{SaWang}
Y.~Sakellaridis and J.~Wang, Intersection complexes and unramified
  {$L$}-factors. \emph{J. Amer. Math. Soc.} \textbf{35} (2022), no.~3,
  799--910. \MR{4433079}

\bibitem{Sarnak}
P.~Sarnak, Comments on {R}obert {L}anglands’ lecture: ``{E}ndoscopy and
  beyond''.

\bibitem{Selberg}
A.~Selberg, Bemerkungen \"{u}ber eine {D}irichletsche {R}eihe, die mit der
  {T}heorie der {M}odulformen nahe verbunden ist. \emph{Arch. Math. Naturvid.}
  \textbf{43} (1940), 47--50. \MR{2626}

\bibitem{Shahidi}
F.~Shahidi, Intertwining operators, {$L$}-functions and representation theory,
  1996, lecture notes of the eleventh {KAIST} mathematics workshop.

\bibitem{Souriau}
J.-M. Souriau, \emph{Structure des syst\`emes dynamiques}. Dunod, Paris, 1970.
  \MR{0260238}

\bibitem{Tate}
J.~T. Tate, Fourier analysis in number fields, and {H}ecke's zeta-functions. In
  \emph{Algebraic {N}umber {T}heory ({P}roc. {I}nstructional {C}onf.,
  {B}righton, 1965)}, pp. 305--347, Thompson, Washington, D.C., 1967.
  \MR{0217026}

\bibitem{vdBS1}
E.~P. van~den Ban and H.~Schlichtkrull, The {P}lancherel decomposition for a
  reductive symmetric space. {I}. {S}pherical functions. \emph{Invent. Math.}
  \textbf{161} (2005), no.~3, 453--566. \MR{2181715 (2006i:43011)}

\bibitem{vdBS2}
E.~P. van~den Ban and H.~Schlichtkrull, The {P}lancherel decomposition for a
  reductive symmetric space. {II}. {R}epresentation theory. \emph{Invent.
  Math.} \textbf{161} (2005), no.~3, 567--628. \MR{2181716 (2006g:22008)}

\bibitem{Venkatesh}
A.~Venkatesh, ``{B}eyond endoscopy'' and special forms on {GL}(2). \emph{J.
  Reine Angew. Math.} \textbf{577} (2004), 23--80. \MR{2108212 (2006b:22016)}

\bibitem{Waldspurger}
J.-L. Waldspurger, Sur les valeurs de certaines fonctions {$L$} automorphes en
  leur centre de sym\'{e}trie. \emph{Compositio Math.} \textbf{54} (1985),
  no.~2, 173--242. \MR{783511}

\bibitem{Waldspurger-Plancherel}
J.-L. Waldspurger, La formule de {P}lancherel pour les groupes {$p$}-adiques
  (d'apr\`es {H}arish-{C}handra). \emph{J. Inst. Math. Jussieu} \textbf{2}
  (2003), no.~2, 235--333. \MR{1989693}

\bibitem{Waldspurger-GP}
J.-L. Waldspurger, La conjecture locale de {G}ross-{P}rasad pour les
  repr\'esentations temp\'er\'ees des groupes sp\'eciaux orthogonaux.
  \emph{Ast\'erisque}  (2012), no. 347, 103--165. \MR{3155345}

\bibitem{Wan}
C.~Wan, On a multiplicity formula for spherical varieties. \emph{J. Eur. Math.
  Soc. (JEMS)} \textbf{24} (2022), no.~10, 3629--3678. \MR{4432908}

\bibitem{Xue}
H.~Xue, Epsilon dichotomy for linear models. \emph{Algebra Number Theory}
  \textbf{15} (2021), no.~1, 173--215. \MR{4226986}

\bibitem{Yamana}
S.~Yamana, L-functions and theta correspondence for classical groups.
  \emph{Invent. Math.} \textbf{196} (2014), no.~3, 651--732. \MR{3211043}

\bibitem{Zhang}
W.~Zhang, Automorphic period and the central value of {R}ankin-{S}elberg
  {L}-function. \emph{J. Amer. Math. Soc.} \textbf{27} (2014), no.~2, 541--612.
  \MR{3164988}

\end{thebibliography}









\end{document}